\documentclass{amsart}
\usepackage{amsfonts}
\usepackage{bbm}
\usepackage{mathrsfs}
\usepackage{color}
\usepackage{amsmath}
\usepackage{bigdelim}
\usepackage{multirow}
\usepackage{amssymb}
\usepackage{indentfirst}
\usepackage{CJK}
\usepackage{hyperref}
\usepackage{fancyhdr}
\usepackage{amscd,amssymb,amsmath,graphicx,verbatim,xy, cases}
\usepackage[TS1,OT1,T1]{fontenc}

\newtheorem{theorem}{Theorem}[section]
\newtheorem{Lemma}[theorem]{Lemma}
\newtheorem{Corollary}[theorem]{Corollary}
\newtheorem{Proposition}[theorem]{Proposition}

\theoremstyle{definition}
\newtheorem{Definition}[theorem]{Definition}
\newtheorem{Example}[theorem]{Example}
\newtheorem{Theorem}[theorem]{Theorem}

\theoremstyle{remark}
\newtheorem{Remark}[theorem]{Remark}
\numberwithin{equation}{section}

\usepackage{bm}
\usepackage{xurl}
\newcommand{\Hq}{\mathbb H}
\newcommand{\Rq}{\mathbb R}
\newcommand{\Cq}{\mathbb C}
\newcommand{\Sp}{\operatorname{Sp}}
\newcommand{\Gr}{\operatorname{Gr}}
\newcommand{\rank}{\operatorname{rank}}
\newcommand{\Ran}{\operatorname{Ran}}
\newcommand{\alg}{\operatorname{alg}_{\Rq}}
\newcommand{\tr}{\operatorname{tr}}
\newcommand{\trr}{\operatorname{tr}_{\Rq}}
\newcommand{\diag}{\operatorname{diag}}
\newcommand{\ii}{\mathbf i}
\newcommand{\jj}{\mathbf j}
\newcommand{\kk}{\mathbf k}
\newcommand{\vct}[1]{\bm{#1}}

\begin{document}
	
	\title[Quaternion Matrix Tuples with Common Invariant Subspaces]{Quaternion Matrix Tuples with Common Invariant Subspaces: Dimension and Nonuniqueness}
	
	%----------Author 1
	\author{Cailing Yao}
	\address{Cailing Yao, School of Medical Informatics, Changchun University of Chinese Medicine, 130117, Changchun, P. R. China}
	\email{yaocl@ccucm.edu.cn}
	%----------Author 2
	\author{Bingzhe Hou}
	\address{Bingzhe Hou, School of Mathematics, Jilin University, 130012, Changchun, P. R. China}
	\email{houbz@jlu.edu.cn}
	%----------Author 3
	\author{Yang Zhang*}\thanks{*Corresponding author.}
	\address{Yang Zhang, School of Medical Informatics, Changchun University of Chinese Medicine, 130117, Changchun, P. R. China}
	\email{zhangyang@ccucm.edu.cn}

	
	
	\begin{abstract}
		Let $M_n(\Hq)$ denote the algebra of $n\times n$ quaternion matrices. We study common invariant right subspaces of tuples in $M_n(\Hq)^r$ and determine exact real semialgebraic dimensions for two associated sets. If $n_1,\ldots,n_s$ are positive integers, $n=n_1+\cdots+n_s$, and $q=\sum_{1\leq a<b\leq s}n_an_b$, the set of tuples preserving at least one nested chain with successive increments $n_1,\ldots,n_s$ has dimension $4rn^2-4(r-1)q$. For $r\geq2$ and $1\leq k<n$, the set of tuples with at least two distinct common $k$-dimensional right invariant subspaces has dimension $4rn^2-4(r-1)(k(n-k)+n-1)$. The proofs use finite right-module coordinates and a common local-rank argument.
	\end{abstract}
	\maketitle

\section{Introduction}

Quaternions were introduced by W. R. Hamilton in 1843 and first published in 1844
\cite{Hamilton1844}.  Their noncommutative multiplication led to a systematic
matrix theory involving eigenvalues, canonical forms, decompositions, and
complex representations \cite{Zhang1997,Rodman2014}.  Quaternion matrices now
provide a unified description of coupled real components in spatial rotation
and attitude control \cite{Shoemake1985,Mayhew2011}, color and multichannel
signal processing \cite{LeBihanSangwine2003}, and quaternion neural networks
\cite{Parcollet2018,Chaudhary2024}. Quaternion-valued methods have also been used in medical image segmentation and multi-channel physiological-signal analysis \cite{C2023,F2025}. These applications motivate the study of
shared lower-dimensional structure in matrix families.

Let $\Hq$ be the division algebra of quaternions, and let $M_n(\Hq)$ denote the
algebra of $n\times n$ quaternion matrices.  For
$A_1,\ldots,A_r\in M_n(\Hq)$, the common invariant subspace problem asks for a
right $\Hq$-subspace $W\subseteq\Hq^n$ such that
\begin{equation}\label{eq:intro-inv}
 A_\ell W\subseteq W,\qquad 1\leq\ell\leq r.
\end{equation}
Here matrices act from the left, whereas scalar multiplication in $W$ is from
the right.  This convention is essential because multiplication in $\Hq$ is
not commutative.
Throughout, the word \emph{reducible} has the following specific meaning: a tuple
is reducible if it has a common invariant
right subspace $W$ with $0<\dim_\Hq W<n$. Invariance of $W^\perp$ is not
required. A common \emph{reducing} subspace is one for which both $W$ and
$W^\perp$ are invariant.

Over $\Cq$, simultaneous-triangularization criteria are associated with
McCoy and Radjavi--Rosenthal \cite{McCoy1936,RadjaviRosenthal2000}.
The common-eigenvector criterion of Shemesh and its finite-family extensions
are given in \cite{Shemesh1984,JamiolkowskiPastuszak2015}. Algebraic and
quantifier-elimination approaches to common invariant subspaces were developed
in \cite{ArapuraPeterson2004,Pastuszak2017}. The compound-matrix criterion
for common invariant subspaces is due to Tsatsomeros \cite{Tsatsomeros2001}.
Incidence constructions appear in the study of Kalman
varieties \cite{OttavianiSturmfels2013}. In particular, U. A. First,
Z. Reichstein, and B. Williams obtained the exact complex dimension
$rn^2-(r-1)k(n-k)$ for the locus of $r$-tuples with a common
invariant subspace of dimension $k$ \cite{FirstReichsteinWilliams2022}.
Related stratifications occur for matrix algebras with involution
\cite{NamTanWilliams2023}. Recent work includes spectral criteria for
reducibility \cite{StessinYang2026} and the algebraic geometry of matrix tuples
with linearly dependent invariant subspaces \cite{Bencze2026}.

The underlying parameter spaces and incidence constructions are classical.
Quaternionic flags admit Bruhat decompositions \cite{MareWillems2009}, while
quiver Grassmannians provide a related framework for subrepresentations over
fields \cite{CerulliIrelliFeiginReineke2012}. The latter is a structural
comparison only: results for Dynkin quivers do not directly apply to arbitrary
loop representations or to right $\Hq$-modules.

Over $\Hq$, D.~\v{Z}.~Djokovi\'c and B.~H.~Smith studied simultaneous
triangularization \cite{DjokovicSmith2008}. Theorem 1.2(iii) in
\cite{Yahaghi2018} classifies irreducible real subalgebras of $M_n(\Hq)$,
up to quaternionic similarity, as $M_n(\Rq)$, $M_n(\Cq)$, or $M_n(\Hq)$.
Thus failure to generate the full real algebra $M_n(\Hq)$ does not by itself
imply reducibility. Burnside-type results for more general division rings are
discussed in \cite{RadjabalipourRosenthalYahaghi2004}. In the complex
representation, an invariant complex subspace corresponds to a right
quaternionic subspace precisely when it is stable under the canonical
conjugate-linear quaternionic structure.

The Shemesh-type criterion below is restricted: it detects
common right eigenvectors whose right eigenvalues commute pairwise, equivalently
whose eigenvalues lie in one complex slice. It is therefore not a general
criterion for reducibility over $\Hq$.

The main results are Theorems~\ref{thm:dimension} and
\ref{thm:nonunique}. For a chain with successive dimension increments
$n_1,\ldots,n_s$, set
\[
q=\sum_{1\leq a<b\leq s}n_an_b.
\]
Theorem~\ref{thm:dimension} gives real semialgebraic dimension
$4rn^2-4(r-1)q$ for the set of $r$-tuples preserving at least one chain of this
type. For $r\geq2$, Proposition~\ref{prop:fixed-h} computes the tuple-image
dimension for each fixed intersection dimension of two distinct common
invariant subspaces of dimension $k$. Taking their finite union yields
$4rn^2-4(r-1)(k(n-k)+n-1)$ in Theorem~\ref{thm:nonunique}.

The dimension proofs use finite right-quaternionic coordinate charts to obtain
incidence upper bounds and an injective differential to attain those bounds
in the tuple images. They do not require all incidence fibers to be finite.
The arbitrary-chain and fixed-intersection calculations are the dimension
contributions of this paper. The complex fixed-dimension theorem and the
$r=2$ complete-chain specialization are established results
\cite{FirstReichsteinWilliams2022,DjokovicSmith2008}; detailed comparisons
appear in Remarks~\ref{rem:complex-chain} and~\ref{rem:complex-nonunique}.

Throughout this paper, dimension means real semialgebraic dimension. In Section~2, we
introduce the notation and preliminary facts. In Section~3, the algebraic
criteria and invariant chains are presented. In Section~4, we show the dimension and nonuniqueness
results.

\section{Preliminaries}
This section introduces the notation and conventions used throughout the
paper. It also presents the basic continuity, rank, and coordinate results
required in the subsequent analysis. See \cite{Rodman2014},
\cite{Zhang1997}, and \cite{HornJohnson2013}.

\subsection{Notation and basic definitions}

We first introduce the algebraic and metric notation used throughout the paper and
identify quaternion matrix spaces with finite-dimensional real coordinate
spaces.

Denote by $n$ the matrix order and $r$ the number of matrices in a
tuple.  Both $n$ and $r$ are positive integers.  Unless another range is displayed, the
index $\ell$ runs from $1$ to $r$. The quaternion division algebra is 
\[
 \Hq=\{a+b\ii+c\jj+d\kk:a,b,c,d\in\Rq\},\ \ii^2=\jj^2=\kk^2=\ii\jj\kk=-1.
\]
For $q=a+b\ii+c\jj+d\kk$, write
$\overline q=a-b\ii-c\jj-d\kk$, $\operatorname{Re}q=a$, and
$|q|=(q\overline q)^{1/2}$. The symbols $\Rq$ and $\Cq$ denote the real and
complex fields, respectively. For $\mathbb F\in\{\Rq,\Cq,\Hq\}$ and positive integers $m,n$,
$M_{m,n}(\mathbb F)$ denotes the set of $m\times n$ matrices over $\mathbb F$,
and $M_n(\mathbb F)=M_{n,n}(\mathbb F)$.  The notation
$\mathbb F^{m\times n}$ has the same meaning.  Thus $M_n(\Hq)^r$ is the
Cartesian product of $r$ copies of $M_n(\Hq)$, whose elements are ordered
$r$-tuples of quaternion matrices.  We write $I_n$ for the $n\times n$
identity matrix and $0_{m\times n}$ for the $m\times n$ zero matrix.  Matrix
sizes will not be suppressed when an identity or zero matrix occurs. If
$A=(a_{ij})\in M_{m,n}(\Hq)$, then $A^*=(\overline{a_{ji}})$
 is its conjugate transpose. The symbols $\Ran A$, $\ker A$, and
$\rank A$ denote its right column space, right null space, and quaternionic
column rank. Thus $A$ has full column rank when its columns are right-linearly
independent. The group of invertible matrices in $M_n(\Hq)$ is denoted by
$GL_n(\Hq)$, and $\diag(A_1,\ldots,A_s)$ denotes a block diagonal matrix.
For a basis matrix $B$, its right column space is
$\Ran B=\{B\vct{c}:\vct{c}\in\Hq^k\}$; changing the chosen basis inside
the same right subspace therefore acts by right multiplication
$B\mapsto BR$ with $R\in GL_k(\Hq)$.

Unless otherwise stated, vectors are column vectors. On the right quaternion vector space $\Hq^n$ we denote by
$\langle \vct{x},\vct{y}\rangle=\vct{x}^*\vct{y}$ the inner product of two vectors $\vct{x}$ and $\vct{y}$, which is right-linear in $\vct{y}$, and
$\|\vct{x}\|_2=\langle\vct{x},\vct{x}\rangle^{1/2}$. The notation
$W^\perp$ means the orthogonal
complement of $W$, while $U\oplus V$ denotes an internal direct sum. The
quaternion dimension of a right subspace is written as $\dim_\Hq W$.
Furthermore,
\[
 \Gr_\Hq(k,n)=\{W\subseteq\Hq^n:W\text{ is a right subspace and }
 \dim_\Hq W=k\}
\]
denotes the quaternionic Grassmannian. A matrix $A\in M_n(\Hq)$ is Hermitian if $A^*=A$,
skew-Hermitian if $A^*=-A$, and unitary if $A^*A=AA^*=I_n$. We write
\[
 \Sp(n)=\{U\in M_n(\Hq):U^*U=UU^*=I_n\}.
\]
A Hermitian idempotent, that is, a matrix $P$ satisfying $P=P^*=P^2$, is the
orthogonal projection onto $\Ran P$.

For a square quaternion matrix define its real trace by
\[
 \trr A=\operatorname{Re}\sum_i a_{ii}.
\]
Only the real trace, rather than the quaternion-valued sum of diagonal
entries, will be used in trace identities. For $X\in M_{m,n}(\Hq)$ and
$Y\in M_{n,m}(\Hq)$, one has $\trr(XY)=\trr(YX)$, since
$\operatorname{Re}(ab)=\operatorname{Re}(ba)$ for $a,b\in\Hq$ and the two
traces are sums over the same pairs of entries. In particular,
$\langle X,Y\rangle_F:=\trr(X^*Y)$ is the real Frobenius inner product.
For a complex matrix, $\tr A$ denotes the ordinary complex trace. For $A=(a_{ij})\in M_{m,n}(\Hq)$ put
\[
 \|A\|_F^2=\sum_{i,j}|a_{ij}|^2.
\]
After writing each quaternion entry in its four real components,
$M_{m,n}(\Hq)$ is a real Euclidean space of dimension $4mn$.
For a tuple $\mathcal A=(A_1,\ldots,A_r)$, the product Frobenius norm is
\[
 \|\mathcal A\|_F^2=\sum_{\ell=1}^r\|A_\ell\|_F^2.
\]

For a subset $S$ of a right vector space over $\mathbb F$,
$\operatorname{span}_{\mathbb F}S$ denotes its right linear span.  For a map
$f$, the symbol $\operatorname{im}f$ denotes its image.  For a complex square
matrix $G$, $\operatorname{spec}(G)$ denotes its complex spectrum, that is, its
set of eigenvalues.  The notation
$\alg(I_n,A_1,\ldots,A_r)$ means the smallest unital real subalgebra of
$M_n(\Hq)$ containing $A_1,\ldots,A_r$.  Equivalently, it consists of all real
linear combinations of finite products of these matrices and $I_n$. We use
$[X,Y]=XY-YX$ to denote the commutator. For a matrix-valued function $R(t)$, the notation
 $R(t)=O(t^2)~ (t\to0)$ means that there exist constants $C>0$ and $\delta>0$, independent of $t$,
such that $\|R(t)\|_F\leq C|t|^2
~\text{whenever } |t|<\delta$. For a differentiable map $F$ between finite-dimensional real coordinate
spaces, $DF(\vct{x})$ denotes its derivative (Jacobian linear map) at
$\vct{x}$.

Finally, a subset $S\subseteq\mathbb R^N$ is called \emph{semialgebraic}
if it can be written as a finite union of sets of the form
\[
\left\{\boldsymbol{x}\in\mathbb R^N:
f_1(\boldsymbol{x})=\cdots=f_p(\boldsymbol{x})=0,\quad
g_1(\boldsymbol{x})>0,\ldots,g_q(\boldsymbol{x})>0
\right\},
\]
where $f_1,\ldots,f_p,g_1,\ldots,g_q$ are real polynomials in $N$
variables. For a nonempty semialgebraic set $S\subseteq\mathbb R^N$, its real
dimension, denoted by $\dim_{\mathbb R}S$, is the largest integer
$d\in\{0,\ldots,N\}$ for which there exists a coordinate projection
$\pi:\mathbb R^N\to\mathbb R^d$ such that $\pi(S)$ contains a nonempty
open subset of $\mathbb R^d$. Here a coordinate projection retains
$d$ selected coordinates. We adopt the convention
$\dim_{\mathbb R}\varnothing=-1$.
See \cite{BasuPollackRoy2006} for the standard theory.
This coordinate-projection definition is equivalent to the standard
semialgebraic dimension defined by maximal-dimensional cells (equivalently, by
the dimensions of coordinate projections). It is defined for arbitrary
semialgebraic sets, including singular ones and it agrees with the manifold
dimension on smooth semialgebraic strata, but
no manifold structure is assumed in the arguments below.

The term locus denotes a set satisfying specified conditions.
For ordered pairs of subspaces $(W,V)$, we use stratum for
 the subset on which the intersection dimension
 $\dim_{\Hq}(W\cap V)$ has a fixed value. A $C^1$ map is a continuously
differentiable map.

\subsection{Continuity, coordinates, and dimension}

In this subsection, we present the continuity, rank, graph-coordinate, and
dimension facts that will be used in the parameter-count arguments.

The following elementary continuity lemma will be used repeatedly.
Its proof is included for completeness.
\begin{Lemma}\label{lem:matrixlimits}
Let $A_j,A,B_j,B\in M_{m,n}(\Hq)$ and
$C_j,C\in M_{n,p}(\Hq)$. If
\[
 \|A_j-A\|_F\longrightarrow0,\qquad
 \|B_j-B\|_F\longrightarrow0,\qquad
 \|C_j-C\|_F\longrightarrow0,
\]
then all matrix-valued limits below are limits in the indicated Frobenius
norms:
\[
 \|(A_j+B_j)-(A+B)\|_F\longrightarrow0,
\]
\[
 \|A_jC_j-AC\|_F\longrightarrow0,
 \qquad
 \|A_j^*-A^*\|_F\longrightarrow0.
\]
If $m=n$, then the real trace
$\trr A=\operatorname{Re}\sum_i a_{ii}$ satisfies
\[
 |\trr A_j-\trr A|\longrightarrow0,
\]
where this is ordinary convergence in $\Rq$.
Consequently, every expression obtained by finitely many additions,
multiplications, adjoints, and real traces passes to the limit termwise.  Its
matrix-valued components converge in Frobenius norm and its real-valued
components converge in the usual absolute-value metric.
\end{Lemma}

\begin{proof}
In the four real coordinates of each quaternion entry, matrix addition,
multiplication, conjugate transpose, and real trace are polynomial maps and
therefore continuous. Frobenius convergence is exactly Euclidean convergence
of these coordinates. The assertions follow, and continuity is preserved
under finite composition.
\end{proof}

\begin{Lemma}[{\cite{Rodman2014}; \cite{Zhang1997}}]\label{lem:graphfacts}
Let $U,V$ be finite-dimensional right quaternionic subspaces. Then
\begin{equation}\label{eq:subspacesum-prelim}
 \dim_\Hq(U+V)=\dim_\Hq U+\dim_\Hq V-\dim_\Hq(U\cap V).
\end{equation}
Furthermore, $\Gr_\Hq(k,N)$ is covered by at most $\binom Nk$ coordinate
regions. In each region a $k$-subspace has the unique form
\begin{equation}\label{eq:graph-prelim}
 W=\Pi^{-1}\Ran\begin{pmatrix}I_k\\ Z\end{pmatrix},
 \ Z\in M_{N-k,k}(\Hq),
\end{equation}
where $\Pi$ is a fixed real row-permutation matrix, namely a real $0$--$1$
matrix whose left action permutes rows. Thus every region has
$k(N-k)$ quaternion coordinates, or $4k(N-k)$ real coordinates.
\end{Lemma}

\begin{proof}
For \eqref{eq:subspacesum-prelim}, extend a basis of $U\cap V$ separately to
bases of $U$ and $V$. The union of the two extensions together with the common
basis is a basis of $U+V$.  The usual coefficient argument is valid over the
division ring $\Hq$ and gives the formula.

For the coordinate statement, let $B\in M_{N,k}(\Hq)$ be a basis matrix of
$W$. Quaternionic Gaussian elimination provides $k$ rows forming an invertible
$k\times k$ submatrix $X$. Move these rows to the top by $\Pi$ and write
$\Pi B=\binom XY$. Right multiplication by $X^{-1}$ leaves the range unchanged
and gives $\binom{I_k}{YX^{-1}}$. Hence $Z=YX^{-1}$. Comparing the top $k$
coordinates proves uniqueness. There are at most $\binom Nk$ possible row sets.
\end{proof}

\begin{Lemma}\label{lem:rankopen}
Let $B_1\in M_{N,p}(\Hq)$ and $B_2\in M_{N,q}(\Hq)$ have full column rank.
Then
\[
 \Ran B_1\cap\Ran B_2=\{\vct{0}\}
 \quad\Longleftrightarrow\quad [\,B_1\ B_2\,]
 \text{ has full column rank}.
\]
If $[\,B_1\ B_2\,]$ has full column rank, then this property,
and hence the triviality of the intersection of the two column spaces, 
persists under all sufficiently small perturbations of $B_1$ and $B_2$.
\end{Lemma}

\begin{proof}
A nonzero vector in the intersection is equivalent to a nontrivial relation
$B_1\vct{x}+B_2\vct{y}=\vct{0}$, which proves the rank criterion. If the
concatenated matrix has
full column rank, some square row submatrix $C$ is invertible. For a sufficiently
 small Frobenius-norm perturbation $\Delta C$, one has
 $\|C^{-1}\Delta C\|_F<1$. The Neumann series
 $\sum_{j=0}^{\infty}(-C^{-1}\Delta C)^j$ is then the inverse of
 $I_{p+q}+C^{-1}\Delta C$. Hence $C+\Delta C$ remains
invertible, and the perturbed concatenation remains full rank.
\end{proof}

\begin{Lemma}[{\cite{BasuPollackRoy2006}}]
\label{lem:semialgebraic-dimension}
Let $X,X_1,\ldots,X_m$ be semialgebraic sets.
\begin{enumerate}
 \item If $f:X\to\Rq^N$ is semialgebraic, then $f(X)$ is semialgebraic and
 $\dim_\Rq f(X)\leq\dim_\Rq X$.
 \item
 $\dim_\Rq\bigl(\bigcup_{j=1}^mX_j\bigr)
 =\max_{1\leq j\leq m}\dim_\Rq X_j$.
 \item If $U\subseteq\Rq^d$ is open and semialgebraic, and
 $f:U\to\Rq^N$ is a semialgebraic $C^1$ map whose derivative is injective at
 one point, then $\dim_\Rq f(U)\geq d$.
\end{enumerate}
The first assertion includes the Tarski--Seidenberg projection theorem, and
the third is the local full-rank criterion used for the lower bounds below.
\end{Lemma}

\begin{proof}
The first two assertions are standard semialgebraic image and dimension
properties, including the Tarski--Seidenberg theorem
\cite{BasuPollackRoy2006}. For (3), let
$\vct{x}_0\in U$ be a point where $Df(\vct{x}_0)$ has rank $d$. Select $d$
output coordinates whose Jacobian minor is invertible, and let $\pi$ retain
those coordinates. The real inverse function theorem implies that
$\pi\circ f$ maps a neighborhood of $\vct{x}_0$ onto an open subset of
$\Rq^d$. The dimension definition therefore gives $\dim_\Rq f(U)\geq d$.
\end{proof}

\section{Algebraic criteria, complex representations, and invariant chains}

In this section, we first present equivalent descriptions of a common
invariant right subspace. We recall the Burnside-type classification of
irreducible real subalgebras, prove a restricted Shemesh-type criterion, and
analyze the structured complex representation. We then identify complete
common invariant chains with simultaneous triangularization.

\subsection{Equivalent formulations and a Burnside-type classification}

We first characterize common invariance in terms of basis matrices, block upper
triangular forms, and Hermitian projections.

As in the introduction, we distinguish between invariant and reducing
subspaces. For a Hermitian idempotent $P$, the quaternionic spectral
theorem gives $\rank P=\trr P$.

The following proposition records standard equivalent formulations of
common invariance. For the underlying quaternionic linear algebra, see \cite{Zhang1997,Rodman2014}. We include a proof for completeness.
\begin{Proposition}\label{prop:four}
Let $1\leq k<n$.  The following are equivalent.
\begin{enumerate}
 \item $\mathcal A$ has a common invariant right subspace of dimension $k$.
 \item There are a full-column-rank $X\in\Hq^{n\times k}$ and matrices
 $B_\ell\in M_k(\Hq)$, uniquely determined by $X$, such that
 \begin{equation}\label{eq:AXXB}
  A_\ell X=XB_\ell\qquad(1\leq\ell\leq r).
 \end{equation}
 \item There exists $U\in\Sp(n)$ such that the matrices have the common block form
 \begin{equation}\label{eq:blockup}
  U^*A_\ell U=
  \begin{pmatrix}
   B_\ell&C_\ell\\
   0_{(n-k)\times k}&D_\ell
  \end{pmatrix},
 \end{equation}
 where $B_\ell$ is the $k\times k$ upper-left block.
 \item There is a $P\in M_n(\Hq)$ such that
 \begin{equation}\label{eq:projection}
  \begin{aligned}
   P^2=P=P^*,\qquad \rank_{\Hq}P=k,
   \quad (I_n-P)A_\ell P=0_{n\times n}\qquad(1\leq\ell\leq r).
  \end{aligned}
 \end{equation}
\end{enumerate}
\end{Proposition}
\begin{proof}
	\emph{$(1)\Longleftrightarrow(2)$.}
	Assume (1), and let the columns of
	$X\in M_{n,k}(\Hq)$ form a right basis of the common invariant
	subspace $W$. For each $\ell$, every column of $A_\ell X$
	belongs to $W$ and therefore has a unique expansion in this basis.
	The corresponding coefficients form a unique matrix
	$B_\ell\in M_k(\Hq)$ satisfying \eqref{eq:AXXB}.
	Conversely, assume (2) and set $W=\Ran X$.
	Since $X$ has full column rank, $\dim_{\Hq}W=k$.
	Moreover, \eqref{eq:AXXB} gives
	\[
	A_\ell W=\Ran(A_\ell X)
	=\Ran(XB_\ell)\subseteq\Ran X=W
	\qquad(1\leq\ell\leq r).
	\]
	Thus (1) holds.
	
	\medskip
	\noindent
	\emph{$(1)\Longleftrightarrow(3)$.}
	Assume (1). Choose an orthonormal basis of $W$ and complete it
	to an orthonormal basis of $\Hq^n$. Write the resulting unitary
	matrix as $U=(U_1\ U_2)$, where $\Ran U_1=W$.
	Since $A_\ell U_1$ has all its columns in $W$ and
	$\Ran U_2=W^\perp$, we have
	\[
	U_2^*A_\ell U_1=0_{(n-k)\times k}
	\qquad(1\leq\ell\leq r).
	\]
	Hence $U^*A_\ell U$ has the form \eqref{eq:blockup}, proving (3).
	Conversely, assume (3), and let $U_1$ consist of the first $k$
	columns of $U$. Multiplying \eqref{eq:blockup} by $U$ on the left
	and comparing the first $k$ columns yields
	\[
	A_\ell U_1=U_1B_\ell
	\qquad(1\leq\ell\leq r).
	\]
	Thus $\Ran U_1$ is a common invariant right subspace of
	dimension $k$, proving (1).
	
	\medskip
	\noindent
	\emph{$(1)\Longleftrightarrow(4)$.}
	Assume (1), and let $P$ be the orthogonal projection onto $W$.
	Then
	\[
	P^2=P=P^*,\qquad
	\rank_{\Hq}P=\trr P=k.
	\]
	Since $A_\ell\Ran P\subseteq\Ran P$ and $P$ acts as the identity
	on $\Ran P$, we obtain
	\[
	PA_\ell P=A_\ell P
	\qquad(1\leq\ell\leq r).
	\]
	Therefore \eqref{eq:projection} holds, proving (4).
	Conversely, assume (4) and set $W=\Ran P$.
	Then $\dim_{\Hq}W=k$. For every $\vct{w}\in W$, we have
	$P\vct{w}=\vct{w}$, and \eqref{eq:projection} implies
	\[
	A_\ell\vct{w}
	=A_\ell P\vct{w}
	=PA_\ell P\vct{w}
	\in\Ran P=W.
	\]
	Hence $W$ is invariant under every $A_\ell$, proving (1).
\end{proof}

\begin{Remark}\label{rem:polynomial-feasibility}
For a Hermitian idempotent $P$, the condition $\rank_{\Hq}P=k$
is equivalent to $\trr P=k$. With this replacement, the conditions in
\eqref{eq:projection} form a finite system of real polynomial equations
in the four real components of the quaternion entries. Hence the existence
of a common invariant right subspace of dimension $k$ is exactly a
real-algebraic feasibility problem. This reformulation gives an exact
algebraic characterization but does not provide an efficient numerical
algorithm for large matrices.
\end{Remark}

Define the unital real matrix algebra
\begin{equation}\label{eq:genalg}
 \mathfrak A=\alg(I_n,A_1,\ldots,A_r).
\end{equation}
A right subspace is invariant under every $A_\ell$ if and only if it is
invariant under every element of $\mathfrak A$.  Since
$\dim_\Rq M_n(\Hq)=4n^2$, this algebra is obtained by a finite linear-span
closure.

The following classification is an existing quaternionic Burnside-type result.
Its necessity is Theorem 1.2(iii) in \cite{Yahaghi2018}.

Two unital real subalgebras $\mathfrak B,\mathfrak C\subseteq M_n(\Hq)$ are
called \emph{quaternionically similar} if there is an $S\in GL_n(\Hq)$ such
that
\begin{equation}\label{eq:algebra-similarity}
 S^{-1}\mathfrak B S
 :=\{S^{-1}BS:B\in\mathfrak B\}=\mathfrak C.
\end{equation}
In this case we write $\mathfrak B\sim_{\Hq}\mathfrak C$.  Thus one quaternion
change of basis simultaneously conjugates every element of the algebra.

We recall the following classification, see Theorem~1.2(iii)
in \cite{Yahaghi2018}.
\begin{Proposition}[{\cite{Yahaghi2018}}]
\label{prop:burnside-classification}
Let $n\geq2$, and let $\mathfrak B\subseteq M_n(\Hq)$ be a unital real
subalgebra. Then $\mathfrak B$ is irreducible on the right quaternionic vector
space $\Hq^n$ if and only if it is quaternionically similar to one of
\begin{equation}\label{eq:three-types}
  M_n(\Rq),\qquad M_n(\Cq),\qquad M_n(\Hq),
\end{equation}
where $\Cq=\Rq+\Rq\ii$ is the fixed complex slice.
\end{Proposition}

\begin{Remark}[Properness and reducibility]
	\label{rem:burnside-classification}
	A proper unital real subalgebra of $M_n(\Hq)$ need not be reducible.
	Indeed, $M_n(\Rq)$ and $M_n(\Cq)$ are proper subalgebras
	that act irreducibly on the right quaternionic vector space $\Hq^n$.
\end{Remark}

\subsection{A Shemesh-type criterion for commuting right eigenvalues}

The next criterion detects a common right eigenvector whose corresponding
right eigenvalues commute pairwise.  The commutativity requirement is exact
and cannot be suppressed.  Its construction is based on the complex
commutator-kernel method of D. Shemesh and its finite-family extensions
\cite{Shemesh1984,JamiolkowskiPastuszak2015}.

Regard $\mathfrak A$ as a finite-dimensional real vector space and choose a
basis
\[
 \mathcal T=(T_1,\ldots,T_d),\qquad
 \mathfrak A=\left\{\sum_{\alpha=1}^d s_\alpha T_\alpha:
 s_\alpha\in\Rq\right\}.
\]
For $[X,Y]=XY-YX$, define
\begin{equation}\label{eq:commutator-kernel}
 \mathcal C(\mathfrak A)
 =\bigcap_{\alpha,\beta=1}^d\ker[T_\alpha,T_\beta]
\end{equation}
and the Hermitian positive semidefinite quaternion matrix
\begin{equation}\label{eq:KT}
 K_{\mathcal T}
 =\sum_{\alpha,\beta=1}^d
 [T_\alpha,T_\beta]^*[T_\alpha,T_\beta].
\end{equation}
The subscript $\mathcal T$ records that the matrix in \eqref{eq:KT} depends
on the selected basis.  Its kernel and singularity are basis-independent.

\begin{Lemma}\label{lem:commutator-kernel}
	The subspace $\mathcal C(\mathfrak A)$ is independent of the choice
	of the real basis $\mathcal T$ of $\mathfrak A$, and
\begin{equation}\label{eq:C-equalities}
 \mathcal C(\mathfrak A)
 =\bigcap_{X,Y\in\mathfrak A}\ker[X,Y]
 =\ker K_{\mathcal T}.
\end{equation}
Moreover, $\mathcal C(\mathfrak A)$ is a right invariant subspace of
$\mathfrak A$.
\end{Lemma}

\begin{proof}
Write
\[
 X=\sum_{\alpha=1}^d s_\alpha T_\alpha,\qquad
 Y=\sum_{\beta=1}^d t_\beta T_\beta,\qquad
 s_\alpha,t_\beta\in\Rq.
\]
For every $\vct{x}\in\Hq^n$,
\begin{equation}\label{eq:commutator-expansion}
 [X,Y]\vct{x}
 =\sum_{\alpha,\beta=1}^d
 s_\alpha t_\beta[T_\alpha,T_\beta]\vct{x}.
\end{equation}
This proves the first equality in \eqref{eq:C-equalities} and its independence
of $\mathcal T$.  Also,
\begin{equation}\label{eq:KT-norm}
 \vct{x}^*K_{\mathcal T}\vct{x}
 =\sum_{\alpha,\beta=1}^d
 \|[T_\alpha,T_\beta]\vct{x}\|_2^2.
\end{equation}
Consequently,
\[
 K_{\mathcal T}\vct{x}=\vct{0}
 \quad\Longleftrightarrow\quad
 [T_\alpha,T_\beta]\vct{x}=\vct{0}
 \quad(1\leq\alpha,\beta\leq d),
\]
Indeed, the forward implication follows from \eqref{eq:KT-norm}, because a
finite sum of nonnegative real numbers is zero only when every summand is
zero.  The reverse implication follows directly from \eqref{eq:KT}.  This
proves the second equality in \eqref{eq:C-equalities}.

Finally, for $R,X,Y\in\mathfrak A$,
\[
 [X,Y]R=[X,YR]-Y[X,R].
\]
If $\vct{x}\in\mathcal C(\mathfrak A)$, then
\[
 [X,Y](R\vct{x})=[X,YR]\vct{x}-Y[X,R]\vct{x}=\vct{0}.
\]
Thus $R\vct{x}\in\mathcal C(\mathfrak A)$.
\end{proof}

The following lemma is an immediate consequence of simultaneous
Schur triangularization for commuting quaternion matrices
\cite{Loring2012}.
\begin{Lemma}
	\label{lem:common-right-eigenvector}
	Let $m,p$ be positive integers, and let
	$B_1,\ldots,B_p\in M_m(\Hq)$ commute pairwise.
	Then there exist
	$\vct{v}\in\Hq^m\setminus\{\vct{0}\}$ and
	$q_1,\ldots,q_p\in\Hq$ such that
	\[
	B_j\vct{v}=\vct{v}q_j
	\qquad(1\leq j\leq p).
	\]
\end{Lemma}

\begin{Theorem}[{Shemesh-type criterion with commuting right eigenvalues}]
\label{thm:shemesh-quaternion}
The following statements are equivalent.
\begin{enumerate}
 \item [(1)] There are $\vct{x}\in\Hq^n\setminus\{\vct{0}\}$ and
 $q_1,\ldots,q_r\in\Hq$ such that
 \begin{equation}\label{eq:commuting-eigenvalues}
  A_\ell\vct{x}=\vct{x}q_\ell\quad(1\leq\ell\leq r),
  \qquad
  q_iq_j=q_jq_i\quad(1\leq i,j\leq r).
 \end{equation}
 \item [(2)] $\mathcal C(\mathfrak A)\ne\{\vct{0}\}$.
 \item [(3)] $K_{\mathcal T}$ is singular (not invertible in $M_n(\Hq)$), or equivalently
 $\ker K_{\mathcal T}\ne\{\vct{0}\}$.
\end{enumerate}
\end{Theorem}
\begin{proof}
	We prove the three implications in turn.
	
	\emph{(1)$\Rightarrow$(2).}
	Assume that there exist
	$\vct{x}\in\Hq^n\setminus\{\vct{0}\}$ and
	$q_1,\ldots,q_r\in\Hq$ satisfying
	\eqref{eq:commuting-eigenvalues}. Since the $q_\ell$ commute pairwise,
	the real algebra
	\[
	\mathcal D=\alg(1,q_1,\ldots,q_r)\subseteq\Hq
	\]
	is commutative. Let
	\[
	W=A_{i_1}\cdots A_{i_m}
	\]
	be a word in the generators, where
	$1\leq i_1,\ldots,i_m\leq r$. Repeatedly applying
	\eqref{eq:commuting-eigenvalues} gives
	\begin{equation}\label{eq:thm31-word-action}
		W\vct{x}
		=
		A_{i_1}\cdots A_{i_m}\vct{x}
		=
		\vct{x}\,q_{i_1}\cdots q_{i_m}.
	\end{equation}
	Every element of $\mathfrak A$ is a real linear combination of such words.
	Consequently, for each $X\in\mathfrak A$ there exists an element
	$s_X\in\mathcal D$ such that
	\begin{equation}\label{eq:thm31-algebra-action}
		X\vct{x}=\vct{x}s_X.
	\end{equation}
	Since $\mathcal D$ is commutative, for all $X,Y\in\mathfrak A$ we obtain
	\begin{equation}\label{eq:thm31-commutator}
		[X,Y]\vct{x}
		=
		X(Y\vct{x})-Y(X\vct{x})
		=
		\vct{x}(s_Xs_Y-s_Ys_X)
		=
		\vct{0}.
	\end{equation}
	Thus $\vct{x}\in\mathcal C(\mathfrak A)$, and hence
	$\mathcal C(\mathfrak A)\neq\{\vct{0}\}$.
	
	\emph{(2)$\Rightarrow$(1).}
	Set $V=\mathcal C(\mathfrak A)$. By Lemma~\ref{lem:commutator-kernel}, $V\neq\{\vct{0}\}$,
	$A_\ell V\subseteq V$ for every $\ell$, and
	\[
	[A_i,A_j]\vct{v}=\vct{0}
	\qquad
	(\vct{v}\in V,\ 1\leq i,j\leq r).
	\]
	Therefore the restrictions $A_\ell|_V$ are well-defined quaternionic
	linear maps on $V$ and commute pairwise,
	\begin{equation}\label{eq:thm31-restricted-commuting}
		[A_i|_V,A_j|_V]=0
		\qquad (1\leq i,j\leq r).
	\end{equation}
	
	Choose a right quaternionic basis of $V$, and let
	$m=\dim_{\Hq}V$. In this basis, the restrictions
	$A_1|_V,\ldots,A_r|_V$ are represented by pairwise commuting
	matrices in $M_m(\Hq)$.
	By Lemma~\ref{lem:common-right-eigenvector}, there exist
	$\vct{x}\in V\setminus\{\vct{0}\}$ and
	$q_1,\ldots,q_r\in\Hq$ such that
	\[
	A_\ell\vct{x}=\vct{x}q_\ell
	\qquad(1\leq\ell\leq r).
	\]
	
	Because $\vct{x}\in V=\mathcal C(\mathfrak A)$, we have
	$[A_i,A_j]\vct{x}=\vct{0}$. Hence
	\begin{equation}\label{eq:thm31-eigenvalue-commutation}
		\begin{aligned}
			\vct{0}
			&=[A_i,A_j]\vct{x} \\
			&=A_i(\vct{x}q_j)-A_j(\vct{x}q_i) \\
			&=\vct{x}q_iq_j-\vct{x}q_jq_i \\
			&=\vct{x}(q_iq_j-q_jq_i).
		\end{aligned}
	\end{equation}
	Choose an index $a$ such that $x_a\neq0$. Taking the $a$-th coordinate
	in \eqref{eq:thm31-eigenvalue-commutation} gives
	\[
	x_a(q_iq_j-q_jq_i)=0.
	\]
	Since every nonzero quaternion is invertible,
	left multiplication by $x_a^{-1}$ yields $q_iq_j=q_jq_i$. Thus condition~(1) holds.
	
	\emph{(2)$\Leftrightarrow$(3).}
	It follows from \eqref{eq:C-equalities} that $\mathcal C(\mathfrak A)=\ker K_{\mathcal T}$. Since $K_{\mathcal T}\in M_n(\Hq)$ is square, it is singular if and only if
	its right kernel is nontrivial. Therefore
	\[
	\mathcal C(\mathfrak A)\neq\{\vct{0}\}
	\quad\Longleftrightarrow\quad
	\ker K_{\mathcal T}\neq\{\vct{0}\}
	\quad\Longleftrightarrow\quad
	K_{\mathcal T}\ \text{is singular}.
	\]
	This proves the equivalence of all three conditions.
\end{proof}

\begin{Remark}[Complex-slice interpretation]
		Quaternions $q_1,\ldots,q_r$ commute pairwise if and only if
	they all belong to a common complex slice $\Rq+\Rq\omega$,
	where $\omega\in\Hq$ satisfies $\omega^2=-1$
	\cite{Zhang1997,Rodman2014}. Thus condition~(1) of Theorem~\ref{thm:shemesh-quaternion}
	is equivalent to the existence of a common right eigenvector
	whose corresponding right eigenvalues belong to a single
	complex slice.

\end{Remark}

\begin{Example}[Necessity of the commutativity condition]
	\label{ex:noncommuting-eigenvalues}
	Let $n=2$ and define
	\[
	A=\ii I_2,\qquad B=\jj I_2.
	\]
	Then the right quaternionic line $W=\operatorname{span}_{\Hq}\{\vct{e}_1\}$ is invariant under both $A$ and $B$, since
	\[
	A\vct{e}_1=\vct{e}_1\ii,\qquad
	B\vct{e}_1=\vct{e}_1\jj.
	\]
	However, $[A,B]
	=(\ii\jj-\jj\ii)I_2
	=2\kk I_2$ is invertible. Consequently,
	\[
	\ker[A,B]=\{\vct{0}\},
	\qquad
	\mathcal C\bigl(\alg(I_2,A,B)\bigr)=\{\vct{0}\}.
	\]
	Thus a common invariant right line may exist even when its corresponding
	right eigenvalues do not commute. Therefore, the commutativity condition
	in Theorem~\ref{thm:shemesh-quaternion} cannot be omitted.
\end{Example}

\subsection{Structured complex representation}
In this subsection, we characterize common invariant right quaternionic
subspaces using the standard complex representation. A common invariant
right subspace of quaternionic dimension $k$ corresponds to a common
invariant complex subspace of dimension $2k$ that is stable under the
conjugate-linear map $\mathcal J$ defined below. For two matrices, we
also give a criterion for the existence of a common invariant right
line in terms of second compound matrices.

Denote $\Cq=\Rq+\Rq\ii$.  Write uniquely $A=C+\jj D$ with
$C,D\in M_n(\Cq)$, and define
\begin{equation}\label{eq:chi}
 \chi(A)=\begin{pmatrix}C&-\overline D\\D&\overline C\end{pmatrix}
 \in M_{2n}(\Cq).
\end{equation}
This injective real-algebra map preserves products and adjoints
\cite{Zhang1997,Rodman2014}.  For
$\vct{x}=\vct{z}+\jj\vct{w}$, put
$\Phi(\vct{x})=(\vct{z}^T,\vct{w}^T)^T$. Here the notation $T$ denotes transpose, and an overline on a
complex vector or matrix denotes entrywise complex conjugation. Right
multiplication by $\jj$ becomes the conjugate-linear map
\begin{equation}\label{eq:J}
 \mathcal J\vct{v}=S\overline{\vct{v}},\qquad
 S=\begin{pmatrix}
  0_{n\times n}&-I_n\\
  I_n&0_{n\times n}
 \end{pmatrix},\qquad \mathcal J^2=-I_{2n}.
\end{equation}
The map $\mathcal J$ is antiunitary. That is, it is conjugate-linear and preserves
the standard complex norm and orthogonality.

The following standard correspondence follows from the complex
	representation of quaternion matrices
	\cite{Zhang1997,Rodman2014}. We include a proof to make the
	dimension correspondence and the equivalent projection conditions
	explicit.
	
	\begin{Theorem}\label{thm:complex}
	
The matrices $A_1,\ldots,A_r\in M_{n}(\Hq)$ have a common invariant right quaternionic
subspace of dimension $k$ if and only if there is a complex subspace
$E\subseteq\Cq^{2n}$ of dimension $2k$ such that
\[
 \chi(A_\ell)E\subseteq E\quad(1\leq\ell\leq r),\qquad
 \mathcal JE=E.
\]
Equivalently, there is a complex Hermitian projection $Q$ satisfying
\begin{equation}\label{eq:complexP}
 Q^2=Q=Q^*,\quad\tr Q=2k,\quad QS=S\overline Q,\quad
 (I_{2n}-Q)\chi(A_\ell)Q=0_{2n\times 2n}.
\end{equation}
\end{Theorem}

\begin{proof}
The identities $c\jj=\jj\overline c$ for $c\in\Cq$ give
\begin{equation}\label{eq:intertwine}
 \Phi(A\vct{x})=\chi(A)\Phi(\vct{x}),\qquad
 \Phi(\vct{x}\jj)=\mathcal J\Phi(\vct{x}).
\end{equation}
Consequently, $\Phi$ identifies right quaternionic subspaces with
$\mathcal J$-stable complex subspaces, preserving matrix invariance. Indeed,
closure under $\Cq$ and under right multiplication by $\jj$ implies closure
under every scalar $a+\jj b$, since
$\vct{x}(a+\jj b)=\vct{x}a+(\vct{x}\jj)b$. If
$\vct{w}_1,\ldots,\vct{w}_k$ is a right $\Hq$-basis of $W$, then
\[
 \Phi(\vct{w}_1),\mathcal J\Phi(\vct{w}_1),\ldots,
 \Phi(\vct{w}_k),\mathcal J\Phi(\vct{w}_k)
\]
is a complex basis of $E=\Phi(W)$: writing coefficients as
$a_i+\jj b_i$, $a_i,b_i\in\Cq$, gives
\[
 \sum_{i=1}^k\bigl(\Phi(\vct{w}_i)a_i+
 \mathcal J\Phi(\vct{w}_i)b_i\bigr)
 =\Phi\!\left(\sum_{i=1}^k\vct{w}_i(a_i+\jj b_i)\right).
\]
Right linear independence proves complex independence, and the same identity
proves spanning. Thus $\dim_\Cq E=2\dim_\Hq W$ in both directions.
For the orthogonal projection $Q$ onto $E$, antiunitarity implies
$\mathcal JE=E$ if and only if the operator identity
$Q\mathcal J=\mathcal JQ$ holds, since $\mathcal J$ then also preserves
$E^\perp$. Here this is an equality of operators on $\mathbb C^{2n}$, not an
ordinary complex-linearity relation. Writing $\mathcal J(v)=S\overline v$,
its matrix form is $QS=S\overline Q$. The remaining equations in
\eqref{eq:complexP} express the dimension and invariance of $\Ran Q$.
\end{proof}

The following example shows that even complex dimension alone does not
imply $\mathcal J$-stability.

\begin{Example}\label{ex:J-constraint}
Let $n=2$ and $A=B=0_{2\times2}$. In the standard complex basis of $\Cq^4$,
set
\[
 E=\operatorname{span}_{\Cq}\{\vct{e}_1,\vct{e}_2\}.
\]
Then $\dim_{\Cq}E=2$, and $E$ is invariant under both $\chi(A)$ and
$\chi(B)$. However,
\[
 \mathcal J\vct{e}_1=\vct{e}_3\notin E,
\]
so $\mathcal JE\ne E$. Hence $E$ is not the image under $\Phi$ of a
right quaternionic line.
\end{Example}
The case $k=1$ of Theorem~\ref{thm:complex} gives the following
equivalent conditions for the existence of a common right eigenvector.

\begin{Corollary}\label{cor:unrestricted-line}
For $A_1,\ldots,A_r\in M_n(\Hq)$, the following statements are equivalent.
\begin{enumerate}
 \item There are $\vct{x}\in\Hq^n\setminus\{\vct{0}\}$ and
 $q_1,\ldots,q_r\in\Hq$ such that
 $A_\ell\vct{x}=\vct{x}q_\ell$ for $1\leq\ell\leq r$.
 \item There is a two-dimensional complex subspace
 $E\subseteq\Cq^{2n}$ satisfying
 \[
  \chi(A_\ell)E\subseteq E\quad(1\leq\ell\leq r),
  \qquad \mathcal JE=E.
 \]
 \item There is a $Q\in M_{2n}(\Cq)$ satisfying
 \[
  Q^2=Q=Q^*,\qquad \tr Q=2,\qquad QS=S\overline Q,
 \]
 and
 \[
  (I_{2n}-Q)\chi(A_\ell)Q=0_{2n\times2n}
  \quad(1\leq\ell\leq r).
 \]
\end{enumerate}
\end{Corollary}	
\begin{proof}
	Condition~(1) is equivalent to the existence of a common invariant
	right line $W=\vct{x}\Hq$. By Theorem~\ref{thm:complex}, such a line
	corresponds to a two-dimensional common invariant complex subspace
	satisfying $\mathcal JE=E$. Thus $(1)\Longleftrightarrow(2)$.
	The equivalence $(2)\Longleftrightarrow(3)$ follows from the
	projection formulation in the same theorem, with $Q$ the
	orthogonal projection onto $E$.
\end{proof}	
For a pair of quaternion matrices, we adapt the complex
compound-matrix criterion of Tsatsomeros~\cite{Tsatsomeros2001}
to the structured complex representation. If $G\in M_{2n}(\Cq)$, its second compound
$C_2(G)$ is the operator induced on the second exterior power
$\bigwedge^2\Cq^{2n}$ by
\[
 C_2(G)(\vct{u}\wedge\vct{v})
 =G\vct{u}\wedge G\vct{v} \qquad
 (\vct{u},\vct{v}\in\Cq^{2n}).
\]

Here $\wedge$ denotes the exterior product. A nonzero vector
$\vct{\xi}\in\bigwedge^2\Cq^{2n}$ is called decomposable if
$\vct{\xi}=\vct{u}\wedge\vct{v}$ for some
$\vct{u},\vct{v}\in\Cq^{2n}$. Since $\vct{\xi}\ne\vct{0}$, the vectors
$\vct{u}$ and $\vct{v}$ are linearly independent and determine the complex plane
\[
E_{\vct{\xi}}
=\operatorname{span}_{\Cq}\{\vct{u},\vct{v}\},
\]
which is independent of the chosen decomposition of $\vct{\xi}$.
By the correspondence in Theorem~\ref{thm:complex},
$E_{\vct{\xi}}$ is the image under $\Phi$ of a right quaternionic
line if and only if $\mathcal J E_{\vct{\xi}}=E_{\vct{\xi}}$.

The second compound construction gives the following criterion for
a common invariant right line of two quaternion matrices.
It requires a decomposable common eigenvector whose associated
complex plane is $\mathcal J$-stable.

For $A,B\in M_n(\Hq)$, choose $\sigma\in\Cq$ such that
$\chi(A)+\sigma I_{2n}$ and $\chi(B)+\sigma I_{2n}$ are
invertible, and set
\[
\widehat A=C_2(\chi(A)+\sigma I_{2n}),\qquad
\widehat B=C_2(\chi(B)+\sigma I_{2n}).
\]
Such a choice is possible because only finitely many values
of $\sigma$ are excluded.

\begin{Corollary}\label{cor:compound-shemesh}
	The matrices $A$ and $B$ have a common invariant right
	quaternionic line if and only if $\widehat A$ and $\widehat B$
	have a nonzero decomposable common eigenvector $\vct{\xi}$
	whose associated complex plane $E_{\vct{\xi}}$ satisfies
	\[
	\mathcal J E_{\vct{\xi}}=E_{\vct{\xi}}.
	\]
\end{Corollary}
\begin{proof}
	Put
	\[
	G=\chi(A)+\sigma I_{2n},\qquad
	H=\chi(B)+\sigma I_{2n}.
	\]
	Both $G$ and $H$ are invertible matrices by the choice of $\sigma$.
	
	Suppose that $A$ and $B$ have a common invariant right line
	$W\subseteq\Hq^n$. Define
	\[
	E:=\Phi(W)=\{\Phi(\vct{x}):\vct{x}\in W\}.
	\]
	By Corollary~\ref{cor:unrestricted-line}, $E$ is a
	$\mathcal J$-stable complex plane and satisfies
	\[
	\chi(A)E\subseteq E,\qquad \chi(B)E\subseteq E.
	\]
	Since
	\[
	G=\chi(A)+\sigma I_{2n},\qquad
	H=\chi(B)+\sigma I_{2n},
	\]
	the scalar shifts preserve $E$, and hence
	\[
	GE\subseteq E,\qquad HE\subseteq E.
	\]
	Because $G$ and $H$ are invertible, these inclusions are equalities:
	\[
	GE=E,\qquad HE=E.
	\]
	Choose a basis $\vct{u},\vct{v}$ of $E$ and set
	\[
	\vct{\xi}=\vct{u}\wedge\vct{v}\ne\vct{0}.
	\]
	Then $\vct{\xi}$ spans $\bigwedge^2E$, which is invariant under
	$C_2(G)$ and $C_2(H)$. Consequently, $\vct{\xi}$ is a nonzero
	decomposable common eigenvector of $\widehat A$ and $\widehat B$,
	and its associated plane satisfies
	\[
	\mathcal J E=E.
	\]
	
	Conversely, let
	$\vct{\xi}=\vct{u}\wedge\vct{v}\ne\vct{0}$
	be a common eigenvector satisfying those conditions, and put
	$E=E_{\vct{\xi}}$. The identity $C_2(XY)=C_2(X)C_2(Y)$ follows from the
	definition. In particular,
	\[
	C_2(G)^{-1}=C_2(G^{-1}),\qquad
	C_2(H)^{-1}=C_2(H^{-1}).
	\]
	
	By the assumption that $\vct{\xi}$ is a common eigenvector of
	$C_2(G)$ and $C_2(H)$, there exist
	$\lambda,\mu\in\Cq$ such that
	\[
	C_2(G)\vct{\xi}=\lambda\vct{\xi},
	\qquad
	C_2(H)\vct{\xi}=\mu\vct{\xi}.
	\]
	Since $C_2(G)$ and $C_2(H)$ are invertible and
	$\vct{\xi}\ne\vct{0}$, both eigenvalues are nonzero. That is, $\lambda\ne0,\ \mu\ne0$. Using $\vct{\xi}=\vct{u}\wedge\vct{v}$ and the definition of the
	second compound matrix, we obtain
	\[
	\begin{aligned}
		G\vct{u}\wedge G\vct{v}
		&=\lambda\,\vct{u}\wedge\vct{v},\\
		H\vct{u}\wedge H\vct{v}
		&=\mu\,\vct{u}\wedge\vct{v}.
	\end{aligned}
	\]

	 Hence
	\[
	\begin{aligned}
		G\vct{u}\wedge G\vct{v}
		&=\lambda\,\vct{u}\wedge\vct{v},\\
		H\vct{u}\wedge H\vct{v}
		&=\mu\,\vct{u}\wedge\vct{v}.
	\end{aligned}
	\]
	
Nonzero scalar multiples of a decomposable vector determine
the same plane. Since $\lambda,\mu\ne0$, the preceding
identities give
\[
\operatorname{span}_{\Cq}\{G\vct{u},G\vct{v}\}
=E
=\operatorname{span}_{\Cq}\{H\vct{u},H\vct{v}\}.
\]
Hence $GE=HE=E$. Since $E$ is a complex subspace, subtracting the scalar shifts
	gives
	\[
	\chi(A)E\subseteq E,\qquad \chi(B)E\subseteq E.
	\]
	Together with $\mathcal JE=E$,
	Corollary~\ref{cor:unrestricted-line} yields a common invariant
	right line for $A$ and $B$.
\end{proof}

For the pair in Example~\ref{ex:J-constraint}, take $\sigma=1$.
Then $\widehat A=\widehat B=I_6$, and
$\vct{\xi}=\vct{e}_1\wedge\vct{e}_2$ is a decomposable common eigenvector.
Its associated plane is the space $E$ in that example, which is not
$\mathcal J$-stable. Thus decomposability alone does not suffice in
Corollary~\ref{cor:compound-shemesh}.

The complex Shemesh construction applied to
$\widehat A$ and $\widehat B$ provides a preliminary criterion. 
A nonzero intersection of their commutator kernels guarantees
the existence of a common complex eigenvector
\cite{Shemesh1984}. In the present quaternionic setting,
this condition is only necessary. The eigenvector must also be
decomposable, and its associated complex plane must satisfy
the $\mathcal J$-stability condition in
Corollary~\ref{cor:compound-shemesh}.

\subsection{Complete invariant chains and simultaneous triangularization}
We now establish the equivalence between the existence of a complete
common invariant chain and simultaneous upper triangularization.
We also record the commuting case. Theorem~2.4(1) in \cite{Loring2012}, stated in terms of structured
complex matrices, yields simultaneous unitary triangularization
of commuting quaternion matrices through the complex representation.
Since unitary similarity is a special case of similarity, the
following proposition is an immediate consequence.

\begin{Proposition}\label{prop:commuting}
If $A_1,\ldots,A_r$ commute pairwise, then one invertible
$S\in GL_n(\Hq)$ makes all $S^{-1}A_\ell S$ upper triangular.
\end{Proposition}
A simultaneous upper triangular form determines a nested sequence
of common invariant right subspaces. We introduce the following
definition to describe this structure.
\begin{Definition}\label{def:completechain}
A nested chain
\[
 \{\vct{0}\}=W_0\subsetneq W_1\subsetneq\cdots\subsetneq W_n=\Hq^n
\]
is called a \emph{complete chain} if
$\dim_\Hq W_a=a$ for $0\leq a\leq n$. It is a complete common invariant
chain for $A_1,\ldots,A_r$ if
$A_\ell W_a\subseteq W_a$ for every $a$ and $\ell$.
\end{Definition}

The following elementary equivalence is obtained by choosing a basis
adapted to the chain.

\begin{Lemma}\label{lem:completechaintriangular}
Let $A_1,\ldots,A_r\in M_n(\Hq)$. Then they possess a complete common invariant
chain if and only if there is an $S\in GL_n(\Hq)$ such that every
$S^{-1}A_\ell S$ is upper triangular.
\end{Lemma}

\begin{proof}
	Suppose that a complete common invariant chain exists.
	Choose a right basis $\vct{v}_1,\ldots,\vct{v}_n$ such that
	\[
	W_a=\operatorname{span}_{\Hq}
	\{\vct{v}_1,\ldots,\vct{v}_a\}
	\qquad(1\leq a\leq n),
	\]
	and put $S=(\vct{v}_1\ \cdots\ \vct{v}_n)$.
	Since $A_\ell\vct{v}_j\in W_j$ for every $\ell$ and $j$,
	each $S^{-1}A_\ell S$ is upper triangular.
	
	Conversely, suppose that $S\in GL_n(\Hq)$ simultaneously
	triangularizes the matrices, and let
	$\vct{v}_1,\ldots,\vct{v}_n$ be its columns.
	Set $W_0=\{\vct{0}\}$ and
	\[
	W_a=\operatorname{span}_{\Hq}
	\{\vct{v}_1,\ldots,\vct{v}_a\}
	\qquad(1\leq a\leq n).
	\]
	Then $\dim_{\Hq}W_a=a$, and upper triangularity gives
	$A_\ell W_a\subseteq W_a$ for every $a$ and $\ell$.
	Thus these subspaces form a complete common invariant chain.
\end{proof}

Proposition~\ref{prop:commuting} and
Lemma~\ref{lem:completechaintriangular} imply that every pairwise
commuting family of quaternion matrices admits a complete common
invariant chain. For further criteria for simultaneous
triangularization, see \cite{DjokovicSmith2008}.

\section{Dimension and nonuniqueness}
In this section, we prove the quaternionic dimension formulas. The dimension formula for complex matrix tuples admitting a common
invariant subspace of prescribed dimension was established in
\cite{FirstReichsteinWilliams2022} and the $r=2$ complete
chain specialization has the quaternionic antecedent
\cite{DjokovicSmith2008}.  We first treat arbitrary prescribed
chain types and then classify pairs of common subspaces by their intersection
dimension before projecting to the tuple space.

Let $\Omega\subseteq\{1,\ldots,n\}\times\{1,\ldots,n\}$
satisfy $(i,i)\in\Omega$ for every $1\leq i\leq n$, and put $c=n^2-|\Omega|$. Define
\[
\begin{aligned}
	\mathcal T_\Omega
	&=\{H\in M_n(\Hq):H_{ij}=0
	\text{ for }(i,j)\notin\Omega\},\\
	\mathcal K_\Omega
	&=\{K\in M_n(\Hq):K_{ij}=0
	\text{ for }(i,j)\in\Omega\}.
\end{aligned}
\]
The entrywise zeros are quaternion scalars. Regarded as a real
vector space, $M_n(\Hq)$ decomposes as
\[
M_n(\Hq)=\mathcal T_\Omega\oplus\mathcal K_\Omega.
\]

Choose pairwise distinct real numbers
$\lambda_1,\ldots,\lambda_n$, and set
$D=\diag(\lambda_1,\ldots,\lambda_n)$.
For $0<\varepsilon<1$, put
\[
U_\varepsilon
=\{K\in\mathcal K_\Omega:\|K\|_F<\varepsilon\}
\times\mathcal T_\Omega^r.
\]
Define $F_\Omega=(F_1,\ldots,F_r)$ on $U_\varepsilon$ by
\begin{equation}\label{eq:separation-map}
	\begin{aligned}
		F_1(K,H_1,\ldots,H_r)
		&=(I_n+K)(D+H_1)(I_n+K)^{-1},\\
		F_\ell(K,H_1,\ldots,H_r)
		&=(I_n+K)H_\ell(I_n+K)^{-1}
		\qquad(2\leq\ell\leq r).
	\end{aligned}
\end{equation}
For $r=1$, only the first component is present.

The following lemma will be used to establish the dimension
lower bounds.

\begin{Lemma}\label{lem:diagonal-separation}
	The map $F_\Omega$ is semialgebraic and continuously
	differentiable. Its derivative at the origin is injective, and
	\[
	\dim_\Rq F_\Omega(U_\varepsilon)
	=4c+4r(n^2-c).
	\]
\end{Lemma}
\begin{proof}
	The real dimensions of $\mathcal K_\Omega$ and
	$\mathcal T_\Omega$ are $4c$ and $4(n^2-c)$, respectively.
	Thus, in real coordinates, $U_\varepsilon$ is an open
	semialgebraic subset of $\Rq^d$, where $d=4c+4r(n^2-c)$.
 For $\|K\|_F<\varepsilon<1$, the Neumann series gives the
	inverse of $I_n+K$. Moreover, for $X\in GL_n(\Hq)$, the adjugate formula gives
	\[
	\chi(X^{-1})
	=\chi(X)^{-1}
	=\frac{\operatorname{adj}(\chi(X))}{\det\chi(X)},
	\]
	where $\operatorname{adj}$ denotes the classical adjugate.
	The real and imaginary parts of the entries of $\chi(X)$
	are real-linear functions of the real coordinates of $X$.
	Consequently, each real coordinate of $X^{-1}$ is a rational
	function of the real coordinates of $X$.
	Thus matrix inversion on $GL_n(\Hq)$ is semialgebraic
	and continuously differentiable. Consequently, $F_\Omega$ is semialgebraic and continuously
	differentiable.
	
	To compute its derivative at the origin, consider
	$F_\Omega(tK,tH_1,\ldots,tH_r)$ for real $t$.
	Using
	\[
	(I_n+tK)^{-1}=I_n-tK+O(t^2)
	\]
	in the defining expressions gives
	\begin{equation}\label{eq:derivative}
		DF_\Omega(\boldsymbol 0)(K,H_1,\ldots,H_r)
		=(H_1+[K,D],H_2,\ldots,H_r),
	\end{equation}
	where $\boldsymbol 0$ denotes the zero parameter tuple.
	
	Suppose that $(K,H_1,\ldots,H_r)$ belongs to the kernel
	of this derivative. Then
	\begin{equation}\label{eq:separation-kernel}
		H_1+[K,D]=0_{n\times n},\qquad
		H_\ell=0_{n\times n}\quad(2\leq\ell\leq r).
	\end{equation}
	Since $D$ is real and diagonal, we obtain
	\[
	[K,D]_{ij}=(\lambda_j-\lambda_i)K_{ij}.
	\]
	Hence $[K,D]\in\mathcal K_\Omega$. As
	$H_1\in\mathcal T_\Omega$ and
	$\mathcal T_\Omega\cap\mathcal K_\Omega
	=\{0_{n\times n}\}$, the first equation in
	\eqref{eq:separation-kernel} implies
	\[
	H_1=0_{n\times n},\qquad [K,D]=0_{n\times n}.
	\]
	For $i\ne j$, the distinctness of the $\lambda_i$ gives
	$K_{ij}=0$. Also, $K_{ii}=0$ because $(i,i)\in\Omega$.
	Thus $K=0_{n\times n}$, and
	\eqref{eq:separation-kernel} shows that every $H_\ell$ vanishes.
	Therefore $DF_\Omega(\boldsymbol 0)$ is injective.
	
	Finally, Lemma~\ref{lem:semialgebraic-dimension}(1) gives
	\[
	\dim_\Rq F_\Omega(U_\varepsilon)\leq d,
	\]
	whereas the injectivity of the derivative and
	Lemma~\ref{lem:semialgebraic-dimension}(3) give the reverse
	inequality. Hence
	\[
	\dim_\Rq F_\Omega(U_\varepsilon)
	=d=4c+4r(n^2-c).
	\]
\end{proof}
In Theorem~\ref{thm:dimension} and
Proposition~\ref{prop:fixed-h}, we apply
Lemma~\ref{lem:diagonal-separation} with $\Omega$
given by the matrix positions allowed by the
corresponding invariance conditions.
\subsection{Prescribed nested chains}

We first construct a finite coordinate cover for nested chains and determine
the block conditions on matrices preserving a fixed chain. These coordinates
yield the exact dimension for every prescribed chain type.

Let $n_1,\ldots,n_s$ be positive integers satisfying
$n_1+\cdots+n_s=n$, and set
\[
 d_0=0,\qquad d_a=n_1+\cdots+n_a\ (1\leq a\leq s),\qquad
 q=\sum_{1\leq a<b\leq s}n_an_b.
\]
We consider chains of right quaternionic subspaces
\begin{equation}\label{eq:chain}
 \{\vct{0}\}=W_0\subset W_1\subset\cdots\subset W_s=\Hq^n,
 \qquad \dim_\Hq W_a=d_a.
\end{equation}
A tuple $\mathcal A=(A_1,\ldots,A_r)$ preserves such a chain if
$A_\ell W_a\subseteq W_a$ for $0\leq a\leq s$ and
$1\leq\ell\leq r$.

Choose an orthonormal basis adapted to the chain, so that $W_a$ is spanned
by its first $d_a$ vectors. A matrix preserves the chain if and only if its
matrix in this basis is block upper triangular with block sizes
$(n_1,\ldots,n_s)$. The strictly lower block part contains
$q$ quaternion entries.  Hence one preserving matrix has
\begin{equation}\label{eq:matrixcount}
 4(n^2-q)
\end{equation}
free real parameters.

We next construct a finite coordinate cover for all chains
of the prescribed type.

Let
\[
 V_a^0=\Hq^{d_a}\oplus\{\vct{0}\}\subseteq\Hq^n
 \qquad(0\leq a\leq s)
\]
be the standard chain, and let
\[
 \mathcal N_-=\{L=(L_{ab})_{a,b=1}^s\in M_n(\Hq):
 L_{ab}=0_{n_a\times n_b}\text{ whenever }a\leq b\}
\]
be the space of strictly block lower triangular matrices. The space $\mathcal N_-$ has $q$ free quaternion entries,
so $\dim_\Rq\mathcal N_-=4q$. Since $L^s=0_{n\times n}$, we have
\begin{equation}\label{eq:lower-inverse-main}
	(I_n+L)^{-1}=\sum_{j=0}^{s-1}(-L)^j.
\end{equation}
Indeed,
\[
(I_n+L)\sum_{j=0}^{s-1}(-L)^j
=\left(\sum_{j=0}^{s-1}(-L)^j\right)(I_n+L)
=I_n-(-L)^s
=I_n.
\]
\begin{Proposition}\label{prop:finite-chain-charts}
	For each real row-permutation matrix $\Pi$, the formulas
	\begin{equation}\label{eq:finite-chain-chart}
		W_a(\Pi,L)=\Pi^{-1}(I_n+L)V_a^0,
		\qquad 0\leq a\leq s,
	\end{equation}
	define an injective map
	$L\mapsto(W_a(\Pi,L))_{a=0}^s$
	from $\mathcal N_-$ to the set of nested chains of right
	quaternionic subspaces with $\dim_\Hq W_a=d_a$.
	Moreover, there is a finite set $\mathfrak S$ of
	row-permutation matrices such that every such chain
	has this form for some $\Pi\in\mathfrak S$
	and $L\in\mathcal N_-$.
\end{Proposition}

\begin{proof}
	The case $s=1$ is immediate, since
	$\mathcal N_-=\{0_{n\times n}\}$ and the only chain is
	$\{\vct{0}\}\subset\Hq^n$. Assume that $s\geq2$. Since $\Pi^{-1}(I_n+L)$ is invertible, it maps the standard
	chain to a nested chain with the prescribed dimensions.
	
	To prove injectivity, fix $\Pi$ and suppose that
	$L,L'\in\mathcal N_-$ represent the same chain. Then
	\[
	C=(I_n+L)^{-1}(I_n+L')
	\]
	preserves every $V_a^0$, so $C$ is block upper triangular.
	On the other hand, $C$ is block lower triangular with
	identity diagonal blocks. Hence $C=I_n$, and therefore
	$L=L'$.
	
	It remains to prove that every chain admits such a
	representation. Let
	$B=[B_1\ \cdots\ B_s]\in GL_n(\Hq)$ be an adapted basis
	matrix, so that its first $d_a$ columns span $W_a$.
	Write $B^{(a)}$ for the matrix consisting of these columns.
	
	We first select nested sets of row indices. By
	Lemma~\ref{lem:graphfacts}, each $B^{(a)}$ contains an
	invertible $d_a\times d_a$ row submatrix. Consequently,
	the rows of $B^{(a)}$ span $\Hq^{1\times d_a}$ as a left
	vector space. Choose a set $I_1$ of $d_1$ row indices
	giving an invertible submatrix of $B^{(1)}$.
	Suppose that $I_a$ has been chosen, where $a<s$.
	The rows of $B^{(a+1)}$ indexed by $I_a$ are
	left-linearly independent. Indeed, a nontrivial left
	linear relation between them would remain valid after
	restriction to the first $d_a$ entries, contradicting
	the invertibility of the previously selected submatrix.
	They can therefore be extended to a basis by adding
	$d_{a+1}-d_a$ rows of $B^{(a+1)}$.
	The resulting index sets satisfy
	\[
	I_1\subset I_2\subset\cdots\subset I_s
	=\{1,\ldots,n\},\qquad |I_a|=d_a.
	\]
	
	Choose a row-permutation matrix $\Pi$ that places the
	rows indexed by $I_a$ in the first $d_a$ positions for
	every $a$. Then the leading principal submatrices of
	$\Pi B$ of orders $d_1,\ldots,d_s$ are invertible.
	
	Block Gaussian elimination now gives the required
	factorization. To justify the successive steps, write
	\[
	\Pi B=
	\begin{pmatrix}
		P_{11}&P_{12}\\
		P_{21}&P_{22}
	\end{pmatrix},
	\qquad P_{11}\in GL_{n_1}(\Hq).
	\]
	Then
	\[
	\Pi B=
	\begin{pmatrix}
		I_{n_1}&0_{n_1\times(n-n_1)}\\
		P_{21}P_{11}^{-1}&I_{n-n_1}
	\end{pmatrix}
	\begin{pmatrix}
		P_{11}&P_{12}\\
		0_{(n-n_1)\times n_1}&Q
	\end{pmatrix},
	\qquad
	Q=P_{22}-P_{21}P_{11}^{-1}P_{12}.
	\]
	For each $2\leq a\leq s$, the leading principal
	submatrix of $Q$ of order $d_a-n_1$ is invertible,
	as follows by applying this factorization to the
	leading principal submatrix of $\Pi B$ of order $d_a$.
	Thus the same argument applies recursively to $Q$.
	It follows that
	\begin{equation}\label{eq:block-elimination-main}
		\Pi B=(I_n+L)R,
	\end{equation}
	where $L\in\mathcal N_-$ and $R$ is invertible and
	block upper triangular. All products occur in the
	displayed order, so the argument is valid over $\Hq$.
Since $RV_a^0=V_a^0$, we obtain
	\[
	W_a=BV_a^0=\Pi^{-1}(I_n+L)V_a^0
	\qquad(0\leq a\leq s).
	\]
	There are only finitely many row-permutation matrices,
	so these representations give the required finite cover.
\end{proof}

The following lemma represents chains by orthogonal
projections and establishes the compactness needed
in the proof of closedness.
\begin{Lemma}\label{lem:compact-chain}
	Assume $s\geq2$ and let
	$\boldsymbol d=(d_1,\ldots,d_{s-1})$. Define
\[
 \mathfrak P_{\boldsymbol d}=
 \left\{(P_1,\ldots,P_{s-1})\in M_n(\Hq)^{s-1}:
 \begin{array}{l}
 P_a^2=P_a=P_a^*,\quad \trr P_a=d_a
 \quad(1\leq a\leq s-1),\\
 P_aP_b=P_a\quad(1\leq a\leq b\leq s-1)
 \end{array}\right\}.
\]
Then $\mathfrak P_{\boldsymbol d}$ is compact in the product Frobenius norm,
and the correspondence $W_a=\Ran P_a$ is a bijection between its elements and
chains of the prescribed dimensions.
\end{Lemma}

\begin{proof}
	For orthogonal projections $P_a$ and $P_b$, the inclusion
	$\Ran P_a\subseteq\Ran P_b$ is equivalent to
	$P_bP_a=P_a$. Taking adjoints shows that this is also
	equivalent to $P_aP_b=P_a$. Moreover,
	\[
	\dim_\Hq\Ran P_a=\rank P_a=\trr P_a=d_a.
	\]
	Since each right subspace has a unique orthogonal
	projection, these observations establish the stated
	bijection.
	
	By Lemma~\ref{lem:matrixlimits}, the defining matrix
	and trace equations are preserved under limits.
	Thus $\mathfrak P_{\boldsymbol d}$ is closed.
	For every
	$\boldsymbol P=(P_1,\ldots,P_{s-1})
	\in\mathfrak P_{\boldsymbol d}$,
	\[
	\|\boldsymbol P\|_F^2
	=\sum_{a=1}^{s-1}\|P_a\|_F^2
	=\sum_{a=1}^{s-1}\trr(P_a^*P_a)
	=\sum_{a=1}^{s-1}d_a.
	\]
	Hence $\mathfrak P_{\boldsymbol d}$ is bounded. As a closed and bounded subset of
	$\Rq^{4n^2(s-1)}$, $\mathfrak P_{\boldsymbol d}$
	is compact.
\end{proof}

The finite coordinate cover in
Proposition~\ref{prop:finite-chain-charts}, together with
Lemma~\ref{lem:semialgebraic-dimension}, yields the
dimension upper bound. A matching lower bound follows
from Lemma~\ref{lem:diagonal-separation}.

Let $\mathfrak F_{\boldsymbol n}$ denote the set of chains
\eqref{eq:chain}, with $\boldsymbol n=(n_1,\ldots,n_s)$. Write
$\boldsymbol W=(W_a)_{a=0}^s$ for such a chain and define
\[
 \widetilde{\mathcal C}_{r,\boldsymbol n}
 =\{(\mathcal A,\boldsymbol W)\in M_n(\Hq)^r\times\mathfrak F_{\boldsymbol n}:
 \ A_\ell W_a\subseteq W_a\text{ for all }a,\ell\}.
\]

We identify each chain $\boldsymbol W$ with the tuple
$(P_1,\ldots,P_{s-1})$, where $P_a$ is the orthogonal
projection onto $W_a$. In these coordinates, the incidence
set $\widetilde{\mathcal C}_{r,\boldsymbol n}$ is defined
by real polynomial equations. The set of matrix tuples
in the following theorem is its image under the
coordinate projection
$\pi(\mathcal A,\boldsymbol W)=\mathcal A$.

\begin{Theorem}\label{thm:dimension}
	Let $\mathcal C_{r,(n_1,\ldots,n_s)}$ denote the set of
	tuples in $M_n(\Hq)^r$ admitting a common invariant chain
	\[
	\{\vct{0}\}=W_0\subset W_1\subset\cdots
	\subset W_s=\Hq^n,
	\qquad
	\dim_\Hq W_a=d_a\quad(0\leq a\leq s).
	\]
	Then $\mathcal C_{r,(n_1,\ldots,n_s)}$ is closed and
	semialgebraic, with
	\begin{equation}\label{eq:dimchain}
		\dim_\Rq\mathcal C_{r,(n_1,\ldots,n_s)}
		=4rn^2-4(r-1)q.
	\end{equation}
	The real codimension in $M_n(\Hq)^r$ is $4(r-1)q$.
\end{Theorem}

\begin{proof}
If $s=1$, then $q=0$ and
$\mathcal C_{r,(n)}=M_n(\Hq)^r$, so all assertions are immediate.
Assume that $s\geq2$. We divide the proof into three steps.

\emph{Step 1: semialgebraicity and the upper bound.}
Let $\mathcal T$ denote the real vector space of block upper triangular
matrices with block sizes $(n_1,\ldots,n_s)$. For each
$\Pi\in\mathfrak S$, set
\[
 S_\Pi(L)=\Pi^{-1}(I_n+L)
\]
and define
\[
 \begin{aligned}
 \Psi_\Pi:\mathcal N_-\times\mathcal T^r&\longrightarrow M_n(\Hq)^r,\\
 (L,T_1,\ldots,T_r)&\longmapsto
 \bigl(S_\Pi(L)T_1S_\Pi(L)^{-1},\ldots,
 S_\Pi(L)T_rS_\Pi(L)^{-1}\bigr).
 \end{aligned}
\]
Since $S_\Pi(L)^{-1}=(I_n+L)^{-1}\Pi$,  the finite expansion of $(I_n+L)^{-1}$ shows that
each real coordinate of $\Psi_\Pi$ is a polynomial
in the real coordinates of $L,T_1,\ldots,T_r$.
Hence $\Psi_\Pi$ is a semialgebraic map. For fixed $L$, a tuple $\mathcal A$ preserves the chain
$(S_\Pi(L)V_a^0)_{a=0}^s$ if and only if
\[
S_\Pi(L)^{-1}A_\ell S_\Pi(L)\in\mathcal T
\qquad(1\leq\ell\leq r).
\]
The finite coordinate cover in
Proposition~\ref{prop:finite-chain-charts} therefore gives
\[
 \mathcal C_{r,(n_1,\ldots,n_s)}
 =\bigcup_{\Pi\in\mathfrak S}\operatorname{im}\Psi_\Pi.
\]
For every $\Pi\in\mathfrak S$, the set
$\Psi_\Pi(\mathcal N_-\times\mathcal T^r)$
is semialgebraic by Lemma~\ref{lem:semialgebraic-dimension}(1).
Since $\mathfrak S$ is finite, their union is semialgebraic.
The domain of each $\Psi_\Pi$ has real dimension
\[
 4q+4r(n^2-q)=4rn^2-4(r-1)q.
\]
The image-dimension and finite-union assertions in
Lemma~\ref{lem:semialgebraic-dimension} give the desired upper bound.

\emph{Step 2: the lower bound.}
Define
\[
\Omega=
\bigcup_{1\leq a\leq b\leq s}
\left\{(i,j):
d_{a-1}<i\leq d_a,\quad
d_{b-1}<j\leq d_b
\right\}.
\] Then
\[
\mathcal T_\Omega=\mathcal T,\qquad
\mathcal K_\Omega=\mathcal N_-,\qquad c=q.
\]
For every $(K,H_1,\ldots,H_r)\in U_\varepsilon$,
the matrices $D+H_1$ and $H_\ell$ $(1\leq\ell\leq r)$
belong to $\mathcal T$. Consequently,
$F_\Omega(K,H_1,\ldots,H_r)$ preserves the chain
\[
\bigl((I_n+K)V_a^0\bigr)_{a=0}^s.
\]
Thus
\[
F_\Omega(U_\varepsilon)
\subseteq\mathcal C_{r,(n_1,\ldots,n_s)}.
\]
By Lemma~\ref{lem:diagonal-separation},
\[
\dim_\Rq\mathcal C_{r,(n_1,\ldots,n_s)}
\geq\dim_\Rq F_\Omega(U_\varepsilon)
=4q+4r(n^2-q).
\]
Together with Step~1, this proves \eqref{eq:dimchain}.

\emph{Step 3: closedness.}
Let
$\mathcal A^{(m)}=(A_1^{(m)},\ldots,A_r^{(m)})
\in\mathcal C_{r,(n_1,\ldots,n_s)}$
converge to $\mathcal A=(A_1,\ldots,A_r)$
in the product Frobenius norm. For each $m$, choose a chain preserved by $\mathcal A^{(m)}$ and let
$P_a^{(m)}$ be the orthogonal projection onto its $a$-th
nontrivial subspace.
Then $(P_1^{(m)},\ldots,P_{s-1}^{(m)})\in\mathfrak P_{\boldsymbol d}$ and
\[
 (I_n-P_a^{(m)})A_\ell^{(m)}P_a^{(m)}=0_{n\times n}
 \qquad(1\leq a\leq s-1,\ 1\leq\ell\leq r).
\]
By Lemma~\ref{lem:compact-chain}, after passing to a subsequence,
these projection tuples converge in the product Frobenius norm to
$(P_1,\ldots,P_{s-1})\in\mathfrak P_{\boldsymbol d}$.
Lemma~\ref{lem:matrixlimits} gives
\[
 (I_n-P_a)A_\ell P_a=0_{n\times n}
 \qquad(1\leq a\leq s-1,\ 1\leq\ell\leq r).
\]
Thus $\mathcal A$ preserves the chain
\[
 \{\vct{0}\}\subset\Ran P_1\subset\cdots
 \subset\Ran P_{s-1}\subset\Hq^n
\]
of the prescribed dimensions. Therefore
$\mathcal A\in\mathcal C_{r,(n_1,\ldots,n_s)}$, proving closedness.
\end{proof}

For $1\leq k<n$, let $\mathcal R_{r,k}$ denote the set of tuples having a
common $k$-dimensional right invariant subspace. Since reducibility in this
paper means the existence of a proper nonzero common invariant right subspace,
define the full reducible locus by
\[
 \mathcal R_r:=\bigcup_{k=1}^{n-1}\mathcal R_{r,k}.
\]
\begin{Corollary}\label{cor:fixedk}
 Let $1\leq k<n$, and let $\mathcal R_{r,k}$ be the set of $r$-tuples with a
 common invariant right quaternionic subspace of dimension $k$.
 Then $\mathcal R_{r,k}$ is closed and semialgebraic, and
\begin{equation}\label{eq:dimR}
 \dim_\Rq\mathcal R_{r,k}
 =4rn^2-4(r-1)k(n-k).
\end{equation}
Moreover, for $n\geq2$, the set $\mathcal R_r$
of reducible $r$-tuples satisfies
\[
\dim_\Rq\mathcal R_r
=4rn^2-4(r-1)(n-1).
\]
For $r\geq2$, $\mathcal R_r$ has positive real codimension
and hence Lebesgue measure zero in $M_n(\Hq)^r$.
\end{Corollary}

\begin{proof}
	For $s=2$, $n_1=k$, and $n_2=n-k$, a chain has
	a single nontrivial subspace of dimension $k$. Thus
	\[
	\mathcal R_{r,k}=\mathcal C_{r,(k,n-k)}.
	\]
	Since $q=k(n-k)$, Theorem~\ref{thm:dimension} shows that
	$\mathcal R_{r,k}$ is closed and semialgebraic
	with dimension given by \eqref{eq:dimR}. The set $\mathcal R_r$ is a finite union of
	semialgebraic sets. By
	Lemma~\ref{lem:semialgebraic-dimension}(2),
	\[
	\dim_\Rq\mathcal R_r
	=\max_{1\leq j<n}\dim_\Rq\mathcal R_{r,j}
	=4rn^2-4(r-1)(n-1).
	\]
The last equality follows from
$\min_{1\leq j<n}j(n-j)=n-1$.
	
	For $r\geq2$, this dimension is strictly less than
	$4rn^2$. Consequently, the semialgebraic set
	$\mathcal R_r$ has Lebesgue measure zero in
	$M_n(\Hq)^r\simeq\Rq^{4rn^2}$.
\end{proof}

\begin{Remark}
	\label{rem:complex-chain}
	Over $\Cq$, the set of $r$-tuples with a common
	$k$-dimensional invariant subspace has complex dimension
	$rn^2-(r-1)k(n-k)$
	\cite{FirstReichsteinWilliams2022}.
	Corollary~\ref{cor:fixedk} gives the corresponding
	real semialgebraic dimension for quaternion matrix tuples.
\end{Remark}

For positive integers $r,n$, define the simultaneously triangularizable locus
\begin{equation}\label{eq:triangularizable-set}
 \mathfrak T_{r,n}:=
 \left\{(A_1,\ldots,A_r)\in M_n(\Hq)^r:
 \begin{array}{l}
 \text{there exists }S\in GL_n(\Hq)\text{ such that}\\
 S^{-1}A_\ell S\text{ is upper triangular for }1\leq\ell\leq r
 \end{array}\right\}.
\end{equation}
By Lemma~\ref{lem:completechaintriangular},
\begin{equation}\label{eq:T-as-chain-locus}
 \mathfrak T_{r,n}=\mathcal C_{r,(1,\ldots,1)}.
\end{equation}

The dimension of the set of simultaneously triangularizable
quaternion matrix pairs was computed in
\cite{DjokovicSmith2008}.
The following corollary gives the formula for arbitrary $r$.
\begin{Corollary}\label{cor:complete-triangularization}
For all positive integers $r,n$,
\begin{equation}\label{eq:dim-triangularizable}
 \dim_\Rq\mathfrak T_{r,n}
 =4rn^2-2(r-1)n(n-1).
\end{equation}
In particular, for $r=2$, the dimension is $2n(3n+1)$.
\end{Corollary}

\begin{proof}
By \eqref{eq:T-as-chain-locus}, it suffices to apply
Theorem~\ref{thm:dimension} to a complete chain. Definition
\ref{def:completechain} gives $s=n$ and $n_a=1$ for every $a$. Hence
$q=\sum_{a<b}1=n(n-1)/2$.  Substitution in \eqref{eq:dimchain} gives
$4rn^2-2(r-1)n(n-1)$. For $r=2$, this simplifies to
$8n^2-2n(n-1)=2n(3n+1)$.
\end{proof}

\subsection{Nonuniqueness of common invariant subspaces}

Let $k$ be an integer with $1\leq k<n$.
We determine the dimension of the set of quaternion
matrix tuples admitting at least two distinct common
invariant right subspaces of dimension $k$.
The proof considers separately the sets of tuples
admitting an ordered pair of such subspaces with
a prescribed intersection dimension and then applies
the finite-union dimension rule.

Let $W,V$ be distinct $k$-dimensional right subspaces and put
$h=\dim_\Hq(W\cap V)$. The dimension formula in
Lemma~\ref{lem:graphfacts} gives
\[
 h_0:=\max\{0,2k-n\}\leq h\leq k-1.
\]
Set
\begin{equation}\label{eq:ch}
 c_h=h(n-h)+2(k-h)(n-k)
 =2k(n-k)-h(n-2k+h).
\end{equation}

Let $\vct{e}_1,\ldots,\vct{e}_n$ be the standard right basis of $\Hq^n$. Put $p=k-h$ and $q_0=n-2k+h$. Define the coordinate subspaces
\[
 \begin{aligned}
 E_0^0&=\operatorname{span}_\Hq\{\vct{e}_1,\ldots,\vct{e}_h\},\\
 E_1^0&=\operatorname{span}_\Hq\{\vct{e}_{h+1},\ldots,\vct{e}_k\},\\
 E_2^0&=\operatorname{span}_\Hq\{\vct{e}_{k+1},\ldots,\vct{e}_{2k-h}\},\\
 E_3^0&=\operatorname{span}_\Hq\{\vct{e}_{2k-h+1},\ldots,\vct{e}_n\},
 \end{aligned}
\]
We adopt the convention
$\operatorname{span}_{\Hq}\varnothing=\{\vct{0}\}$. The
\emph{standard pair with intersection dimension $h$} is
\begin{equation}\label{eq:standard-pair}
 W_0=E_0^0\oplus E_1^0,
 \qquad V_0=E_0^0\oplus E_2^0.
\end{equation}
It satisfies
$\dim_\Hq W_0=\dim_\Hq V_0=k$ and
$\dim_\Hq(W_0\cap V_0)=h$.

Let $h$ be an integer with $h_0\leq h\leq k-1$.
Denote by $\mathcal G_{n,k,h}$ the set of ordered pairs
$(W,V)$ of $k$-dimensional right subspaces of $\Hq^n$
satisfying $\dim_\Hq(W\cap V)=h$.
\begin{Lemma}\label{lem:pair-coordinates}
	The set $\mathcal G_{n,k,h}$ is semialgebraic and has
	real dimension $4c_h$.
	There exist finitely many open semialgebraic sets
	$U_\nu\subseteq\Rq^{4c_h}$ and semialgebraic $C^1$ maps
	$S_\nu:U_\nu\to GL_n(\Hq)$ such that the maps
	\[
	\begin{aligned}
		f_\nu:U_\nu&\longrightarrow\mathcal G_{n,k,h},\\
		\vct{u}&\longmapsto
		\bigl(S_\nu(\vct{u})W_0,\,
		S_\nu(\vct{u})V_0\bigr)
	\end{aligned}
	\]
	are injective and their images cover $\mathcal G_{n,k,h}$.
\end{Lemma}

\begin{proof}
	Put $N=n-h$ and $p=k-h$. For $(W,V)\in\mathcal G_{n,k,h}$,
	set
	\[
	E_0=W\cap V,\qquad
	E_1=W\cap E_0^\perp,\qquad
	E_2=V\cap E_0^\perp.
	\]
	Then
	\[
	W=E_0\oplus E_1,\qquad V=E_0\oplus E_2,
	\qquad E_1\cap E_2=\{\vct{0}\}.
	\]
	These three subspaces are uniquely determined by $(W,V)$.
	
	Suppose first that $h>0$. By
	Lemma~\ref{lem:graphfacts}, finitely many coordinate
	regions suffice to represent $E_0$. In one such region,
	\[
	E_0=\Pi^{-1}\Ran
	\begin{pmatrix}I_h\\ Z\end{pmatrix},
	\qquad Z\in M_{N,h}(\Hq),
	\]
	where $\Pi$ is a fixed row-permutation matrix. Define
	\begin{equation}\label{eq:E0-trivialization}
		G(Z)=
		\begin{pmatrix}
			I_h&-Z^*\\
			Z&I_N
		\end{pmatrix}.
	\end{equation}
	Since
	\[
	G(Z)^*G(Z)
	=\diag(I_h+Z^*Z,I_N+ZZ^*),
	\]
	and both diagonal blocks are positive definite,
	$G(Z)$ is invertible. The first $h$ columns of
	$\Pi^{-1}G(Z)$ span $E_0$, and the remaining $N$ columns
	span $E_0^\perp$. The latter columns therefore identify
	$\Hq^N$ with $E_0^\perp$.
	
	When $h=0$, the subspace $E_0$ is $\{\vct{0}\}$.
	There are no parameters $Z$, and the identification
	of $\Hq^N$ with $E_0^\perp$ is the identity.
	
	Using this identification, apply the finite graph
	coordinates of Lemma~\ref{lem:graphfacts} to $E_1$ and
	$E_2$. Let $B_1,B_2\in M_{N,p}(\Hq)$ be their graph
	basis matrices in $\Hq^N$. The combined real coordinates
	form a vector $\vct{u}\in\Rq^d$, where
	\[
	d=4hN+8p(N-p)=4c_h.
	\]
	By Lemma~\ref{lem:rankopen}, the condition
	$E_1\cap E_2=\{\vct{0}\}$ is equivalent to
	\[
	B=[\,B_1\ B_2\,]
	\]
	having full column rank.
	
	Refine each coordinate region by choosing a set $J$
	of $2p$ row indices and requiring the corresponding
	square submatrix $B_J$ to be invertible.
	These conditions cover all full-column-rank matrices $B$. For each fixed $J$, the parameter values for which
	$B_J$ is invertible form an open semialgebraic set,
	defined by
	\[
	|\det\chi(B_J)|^2>0.
	\]
Thus the resulting parameter domains
	$U_\nu\subseteq\Rq^d$ are open and semialgebraic.
	
	For each fixed $J$, complete $B$ by the standard
	coordinate columns indexed by the complement of $J$,
	and denote the resulting square matrix by $R$.
	If $q_0=N-2p>0$, a row permutation makes $R$ block
	lower triangular with diagonal blocks $B_J$ and
	$I_{q_0}$. Hence $R$ is invertible.
	If $q_0=0$, take $R=B$, which is already invertible.
	
	On each parameter domain, define
	\[
	S_\nu(\vct{u})=
	\begin{cases}
		\Pi^{-1}G(Z)\diag(I_h,R),& h>0,\\
		R,& h=0.
	\end{cases}
	\]
	Then
	\[
	S_\nu(\vct{u})W_0=W,\qquad
	S_\nu(\vct{u})V_0=V.
	\]
	Each real coordinate of $S_\nu(\vct{u})$ is a polynomial
	in the components of $\vct{u}$. Hence $S_\nu$ is
	a semialgebraic $C^1$ map.
		
Each map $f_\nu$ is injective. Indeed,
	$(W,V)$ uniquely determines $E_0,E_1,E_2$, and the
	graph coordinates of these subspaces are unique
	within each fixed coordinate region. Only finitely
	many graph regions and row sets $J$ occur, so the
	images of these maps cover $\mathcal G_{n,k,h}$.
	
Each map $f_\nu$ is semialgebraic when subspaces
are represented by their orthogonal projections. To see this, recall that for a full-column-rank
	matrix $X$, the orthogonal projection onto $\Ran X$ is
	\[
	X(X^*X)^{-1}X^*.
	\]	
	The assertion follows from the matrix-inversion
	argument in Lemma~\ref{lem:diagonal-separation}.
	Each injective semialgebraic pair map has a
	semialgebraic inverse on its image. Applying Lemma~\ref{lem:semialgebraic-dimension}(1)
	to $f_\nu$ and its inverse
	$f_\nu^{-1}:f_\nu(U_\nu)\to U_\nu$ gives
	\[
	\dim_\Rq f_\nu(U_\nu)
	=\dim_\Rq U_\nu
	=4c_h
	\]
	for every nonempty $U_\nu$.
	
	Since $(W_0,V_0)\in\mathcal G_{n,k,h}$, at least one
	parameter domain is nonempty. The finite-union
	property in Lemma~\ref{lem:semialgebraic-dimension}(2)
	now gives
	\[
	\dim_\Rq\mathcal G_{n,k,h}=d=4c_h.
	\]
	The same finite cover proves semialgebraicity.
\end{proof}

Choose a decomposition
\begin{equation}\label{eq:pair-decomposition}
 \Hq^n=E_0\oplus E_1\oplus E_2\oplus E_3,
 \qquad W=E_0\oplus E_1,\quad V=E_0\oplus E_2,
\end{equation}
with block dimensions
\[
 (h,p,p,q_0),\qquad p=k-h,\qquad q_0=n-2k+h.
\]

\begin{Lemma}\label{lem:pair-stabilizer}
With respect to the ordered decomposition
$\Hq^n=E_0\oplus E_1\oplus E_2\oplus E_3$,
a matrix $X\in M_n(\Hq)$ preserves both $W$ and $V$
if and only if it has the block form	
\begin{equation}\label{eq:pair-block}
 X=\begin{pmatrix}
 X_{00}&X_{01}&X_{02}&X_{03}\\
 0_{p\times h}&X_{11}&0_{p\times p}&X_{13}\\
 0_{p\times h}&0_{p\times p}&X_{22}&X_{23}\\
 0_{q_0\times h}&0_{q_0\times p}&0_{q_0\times p}&X_{33}
 \end{pmatrix}.
\end{equation}
Block rows and columns corresponding to zero-dimensional
subspaces are omitted. The real vector space of such matrices has dimension $4(n^2-c_h)$.
\end{Lemma}
\begin{proof}
	Suppose that $XW\subseteq W$ and $XV\subseteq V$.
	Since $E_0=W\cap V$, we obtain
	\[
	XE_0\subseteq E_0,\qquad
	XE_1\subseteq E_0\oplus E_1,\qquad
	XE_2\subseteq E_0\oplus E_2.
	\]
	These inclusions give exactly the zero blocks in the
	stated matrix form. 
	
	Conversely, the block form in (\ref{eq:pair-block}) gives
	\[
	\begin{aligned}
		XW&=X(E_0\oplus E_1)\subseteq E_0\oplus E_1=W,\\
		XV&=X(E_0\oplus E_2)\subseteq E_0\oplus E_2=V.
	\end{aligned}
	\]
	This proves the equivalence.
	
	The first block column has $h(2p+q_0)$ entries
	required to vanish, and each of the next two block
	columns has $p(p+q_0)$ such entries. These positions
	are disjoint, and their total number is
	\[
	h(2p+q_0)+2p(p+q_0)
	=h(n-h)+2(k-h)(n-k)
	=c_h.
	\]
	The remaining $n^2-c_h$ quaternion entries vary
	independently. Since each contributes four real
	coordinates, the resulting real vector space has
	dimension $4(n^2-c_h)$.
\end{proof}

To distinguish the subspace-pair parameters from their tuple image, define
the incidence set
\begin{equation}\label{eq:pair-incidence}
 \widetilde{\mathcal M}_{r,k,h}
 =\left\{(\mathcal A,W,V):
 \begin{array}{l}
 \mathcal A\in M_n(\Hq)^r,\quad (W,V)\in\mathcal G_{n,k,h},\\
 A_\ell W\subseteq W,\quad A_\ell V\subseteq V\quad(1\leq\ell\leq r)
 \end{array}\right\}.
\end{equation}
Let $\pi(\mathcal A,W,V)=\mathcal A$ and set
\begin{equation}\label{eq:pair-image}
 \mathcal M_{r,k,h}=\pi(\widetilde{\mathcal M}_{r,k,h}).
\end{equation}
The set $\mathcal M_{r,k,h}$ consists of matrix tuples
admitting two distinct common invariant right subspaces
$W,V$ with
$\dim_\Hq W=\dim_\Hq V=k$ and
$\dim_\Hq(W\cap V)=h$. A tuple can belong to more
than one of the sets $\mathcal M_{r,k,h}$.

\begin{Proposition}\label{prop:fixed-h}
Let $r\geq2$, $1\leq k<n$, and put
$h_0=\max\{0,2k-n\}$. For each integer
$h_0\leq h\leq k-1$, the image set
$\mathcal M_{r,k,h}=\pi(\widetilde{\mathcal M}_{r,k,h})$ is semialgebraic and
\begin{equation}\label{eq:fixed-h-dimension}
 \dim_\Rq\mathcal M_{r,k,h}=4rn^2-4(r-1)c_h.
\end{equation}
\end{Proposition}

\begin{proof}
Let $\mathcal T_h$ be the real vector space of matrices
preserving both $W_0$ and $V_0$.
By Lemma~\ref{lem:pair-stabilizer},
\[
\dim_\Rq\mathcal T_h=4(n^2-c_h).
\]	
 Let
$U_\nu\subseteq\Rq^{4c_h}$ be one of the finitely many semialgebraic
coordinate regions in Lemma~\ref{lem:pair-coordinates}, and let
$S_\nu(\vct{u})$ be its adapted invertible matrix. Define
\[
 \begin{aligned}
 \Psi_\nu:U_\nu\times\mathcal T_h^r&\longrightarrow M_n(\Hq)^r,\\
 (\vct{u},T_1,\ldots,T_r)&\longmapsto
 \bigl(S_\nu(\vct{u})T_1S_\nu(\vct{u})^{-1},\ldots,
 S_\nu(\vct{u})T_rS_\nu(\vct{u})^{-1}\bigr).
 \end{aligned}
\]
The inversion fact proved in Lemma~\ref{lem:diagonal-separation} makes
$\Psi_\nu$ semialgebraic. Define
\[
\begin{aligned}
	\widetilde\Psi_\nu:
	U_\nu\times\mathcal T_h^r
	&\longrightarrow\widetilde{\mathcal M}_{r,k,h},\\
	(\vct{u},T_1,\ldots,T_r)
	&\longmapsto
	\bigl(
	\Psi_\nu(\vct{u},T_1,\ldots,T_r),\,
	S_\nu(\vct{u})W_0,\,
	S_\nu(\vct{u})V_0
	\bigr).
\end{aligned}
\]
The image of $\widetilde\Psi_\nu$ consists precisely of
the triples $(\mathcal A,W,V)$ in
$\widetilde{\mathcal M}_{r,k,h}$ for which
$(W,V)\in f_\nu(U_\nu)$.

Since $\Psi_\nu$ and $f_\nu$ are semialgebraic,
$\widetilde\Psi_\nu$ is semialgebraic.
The parameters are uniquely recovered from a triple
in its image by
\[
\vct{u}=f_\nu^{-1}(W,V),\qquad
T_\ell=S_\nu(\vct{u})^{-1}A_\ell S_\nu(\vct{u})
\quad(1\leq\ell\leq r).
\]
These formulas show that
$\widetilde\Psi_\nu^{-1}$ is semialgebraic on
$\widetilde\Psi_\nu(U_\nu\times\mathcal T_h^r)$.
Thus every nonempty image of $\widetilde\Psi_\nu$
has the same real dimension as its parameter domain.
Since these images form a finite cover of
$\widetilde{\mathcal M}_{r,k,h}$,
Lemma~\ref{lem:semialgebraic-dimension} gives
\[
\dim_\Rq\widetilde{\mathcal M}_{r,k,h}
=4c_h+4r(n^2-c_h).
\]
For each $\nu$, the set
$\Psi_\nu(U_\nu\times\mathcal T_h^r)$ is semialgebraic
by Lemma~\ref{lem:semialgebraic-dimension}(1).
Since
\[
\mathcal M_{r,k,h}
=\bigcup_\nu\Psi_\nu(U_\nu\times\mathcal T_h^r)
\]
and the union is finite, $\mathcal M_{r,k,h}$
is semialgebraic. The image-dimension and finite-union properties in
Lemma~\ref{lem:semialgebraic-dimension} give
\[
 \dim_\Rq\mathcal M_{r,k,h}
 \leq4c_h+4r(n^2-c_h).
\]

For the reverse inequality, use the standard pair \eqref{eq:standard-pair} and
the coordinate decomposition
$\Hq^n=E_0^0\oplus E_1^0\oplus E_2^0\oplus E_3^0$ defined above. Let $\mathcal K_h$ be the real vector space of quaternion
matrices with arbitrary entries in the blocks required
to vanish in Lemma~\ref{lem:pair-stabilizer}
and zeros in all remaining blocks. The $c_h$ quaternion entries in these blocks vary
independently, so
\[
\dim_\Rq\mathcal K_h=4c_h.
\]
Let $\Omega_h$ be the set of scalar matrix positions
contained in the blocks that may be nonzero in the block
pattern defining $\mathcal T_h$. Then
\[
\mathcal T_{\Omega_h}=\mathcal T_h,\qquad
n^2-|\Omega_h|=c_h,\qquad
\mathcal K_{\Omega_h}=\mathcal K_h.
\]
Applying Lemma~\ref{lem:diagonal-separation} with
$\Omega=\Omega_h$ gives $c=c_h$ and
$\mathcal K_\Omega=\mathcal K_h$. Fix $\varepsilon\in(0,1)$ and define the parameter domain
\[
U_{\varepsilon,h}
=
\{K\in\mathcal K_h:\|K\|_F<\varepsilon\}
\times\mathcal T_h^r.
\]
After identifying
$\mathcal K_h\times\mathcal T_h^r$
with its real coordinate space,
$U_{\varepsilon,h}$ is an open semialgebraic set.
For every
$(K,H_1,\ldots,H_r)\in U_{\varepsilon,h}$,
the matrix
\[
S(K)=I_n+K
\]
is invertible. By Lemma~\ref{lem:diagonal-separation}, applied with
$\Omega=\Omega_h$,
\[
\dim_\Rq F_{\Omega_h}(U_{\varepsilon,h})
=4c_h+4r(n^2-c_h).
\]
For every
$(K,H_1,\ldots,H_r)\in U_{\varepsilon,h}$,
the matrices $D+H_1,H_2,\ldots,H_r$ preserve
the standard pair $(W_0,V_0)$. Hence
$F_{\Omega_h}(K,H_1,\ldots,H_r)$ preserves
\[
\bigl((I_n+K)W_0,\,(I_n+K)V_0\bigr).
\]
Since $I_n+K$ is invertible,
\[
\dim_\Hq\bigl((I_n+K)W_0\cap(I_n+K)V_0\bigr)=h.
\]
Therefore
\[
F_{\Omega_h}(U_{\varepsilon,h})
\subseteq\mathcal M_{r,k,h},
\]
and consequently
\[
\dim_\Rq\mathcal M_{r,k,h}
\geq
\dim_\Rq F_{\Omega_h}(U_{\varepsilon,h})
=
4c_h+4r(n^2-c_h),
\]
which is the required reverse inequality. 
\end{proof}

\begin{Theorem}\label{thm:nonunique}
Let $r\geq2$, $1\leq k<n$, and set
$h_0=\max\{0,2k-n\}$. Let $\mathcal M_{r,k}$ consist of tuples having at least
two distinct common invariant right quaternionic subspaces of dimension $k$.
Then
\begin{equation}\label{eq:dimmultiple}
 \dim_\Rq\mathcal M_{r,k}
 =4rn^2-4(r-1)(k(n-k)+n-1).
\end{equation}
Moreover,
\[
\dim_\Rq\mathcal R_{r,k}
-
\dim_\Rq\mathcal M_{r,k}
=4(r-1)(n-1).
\]
\end{Theorem}

\begin{proof}
By definition,
\[
 \mathcal M_{r,k}=\bigcup_{h=h_0}^{k-1}\mathcal M_{r,k,h}.
\]
The index range is nonempty in all cases. If $2k\leq n$, then
$h_0=0$ and $h=0,\ldots,k-1$. If $2k>n$, then
$h_0=2k-n$ and $h=2k-n,\ldots,k-1$. In particular, for $k=1$ the
only value is $h=0$, whereas for $k=n-1$ the only value is $h=n-2$. The index set $\{h_0,\ldots,k-1\}$ is finite.
Hence $\mathcal M_{r,k}$ is a finite union of the sets
$\mathcal M_{r,k,h}$. By
Lemma~\ref{lem:semialgebraic-dimension}(2) and
Proposition~\ref{prop:fixed-h}, its dimension is determined
by the smallest value of $c_h$. Direct rearrangement gives
\[
 c_h-k(n-k)=(k-h)(n-k+h).
\]
Put $s=k-h$.  The range of $h$ gives
$1\leq s\leq k$ when $2k\leq n$, and
$1\leq s\leq n-k$ when $2k>n$. Thus $s=1$ belongs to the admissible range in both cases. Since $s\leq n-1$ and
\[
 s(n-s)-(n-1)=(s-1)(n-s-1)\geq0,
\]
 the minimum of $s(n-s)$ on this range is $n-1$. The value $s=1$ attains equality. If $s=n-1$ also belongs
 to the relevant range, it gives the same value.
 
Consequently, the minimum of $c_h$ is
$k(n-k)+n-1$, which gives the largest semialgebraic dimension. Taking the union of the finitely many sets
$\mathcal M_{r,k,h}$ preserves this maximum and proves
\eqref{eq:dimmultiple}. Combining this formula with the
dimension formula for $\mathcal R_{r,k}$ in
Corollary~\ref{cor:fixedk} gives
\[
\dim_\Rq\mathcal R_{r,k}
-
\dim_\Rq\mathcal M_{r,k}
=4(r-1)(n-1).
\]
This proves the remaining assertion.

\end{proof}

As a direct consequence, for $r\geq2$ the semialgebraic set
\[
 \mathcal U_{r,k}:=\mathcal R_{r,k}\setminus\mathcal M_{r,k}
\]
of tuples having exactly one common $k$-dimensional right invariant subspace
satisfies
\[
 \dim_\Rq\mathcal U_{r,k}=\dim_\Rq\mathcal R_{r,k}.
\]
Both $\mathcal R_{r,k}$ and $\mathcal M_{r,k}$ are semialgebraic,
and the preceding dimension formulas give
\[
\dim_\Rq\mathcal M_{r,k}
<
\dim_\Rq\mathcal R_{r,k}.
\]
Since
\[
\mathcal R_{r,k}
=\mathcal U_{r,k}\cup\mathcal M_{r,k},
\]
Lemma~\ref{lem:semialgebraic-dimension}(2) yields
\[
\dim_\Rq\mathcal R_{r,k}
=
\max\bigl\{
\dim_\Rq\mathcal U_{r,k},
\dim_\Rq\mathcal M_{r,k}
\bigr\}
=
\dim_\Rq\mathcal U_{r,k}.
\]

The set $\mathcal M_{r,k}$ consists exactly of the
tuples that have at least two distinct common invariant
right subspaces of dimension $k$. The dimension formula
is a global statement about semialgebraic dimension and
does not determine the local behavior near any particular
tuple. In particular, it does not imply that the set of
tuples with more than one common $k$-dimensional invariant
subspace is the same as the set of tuples having infinitely
many such subspaces. Some lower-dimensional subsets of $\mathcal M_{r,k}$
may nevertheless contain tuples that preserve infinitely
many common invariant right subspaces of dimension $k$.

\begin{Remark}[Complex counterpart]
	\label{rem:complex-nonunique}
	For comparison, the same incidence calculation over $\Cq$
	gives
	\[
	\dim_{\Cq}\mathcal M^{\Cq}_{r,k}
	=
	rn^2-(r-1)\bigl(k(n-k)+n-1\bigr),
	\]
	where $\mathcal M^{\Cq}_{r,k}$ denotes the set of complex
	$r$-tuples having at least two distinct common invariant
	subspaces of dimension $k$. The calculation above is the
	corresponding real semialgebraic dimension formula for
	quaternion matrix tuples.
\end{Remark}

\section*{Declarations}

\noindent \textbf{Ethics approval}

\noindent Not applicable.

\noindent \textbf{Competing interests}

\noindent The author declares that there is no conflict of interest or competing interest.

\noindent \textbf{Funding}

\noindent The third author was supported by Jilin Provincial Higher Education Teaching Reform Research Project (Grant No. JLJY202731893340).

\noindent \textbf{Authors' contributions}

\noindent All authors contributed equally to this work.

\noindent \textbf{Availability of data and materials}

\noindent Data sharing is not applicable to this article as no data sets were generated or analyzed during the current study.

\end{document}